\documentclass[bj,preprint]{imsart}

\usepackage{graphicx}
\usepackage[space]{grffile}
\usepackage{latexsym}
\usepackage{textcomp}
\usepackage{longtable}
\usepackage{tabulary}
\usepackage{booktabs,array,multirow}
\usepackage{amsfonts,amsmath,amssymb}
\usepackage[dvipsnames]{xcolor}

\usepackage{hyperref}  % hyperlinks
\hypersetup{
    colorlinks=true,
    linkcolor=blue,
    filecolor=magenta,      
    urlcolor=brown,
    citecolor = red,
    pdftitle={Overleaf Example},
    pdfpagemode=FullScreen,
    }
\usepackage{url}  
\usepackage{etoolbox}
\usepackage{graphicx,psfrag,epsf}
\usepackage{enumerate}
\usepackage{enumitem}
\usepackage{graphics}
\usepackage{amsthm}
\usepackage{multirow}
\usepackage{mathtools}

\DeclarePairedDelimiter\floor{\lfloor}{\rfloor}
\newtheorem{theorem}{Theorem}
\newtheorem{proposition}{Proposition} %% remove [theorem] for Unique numbering
\theoremstyle{plain}
\newtheorem{corro}{Corollary}

\newenvironment{example}[1][Example]{\begin{trivlist}
\item[\hskip \labelsep {\bfseries #1}]}{\end{trivlist}}

\def\Cum{\mbox{Cum}}

\usepackage{mac-09-ams-03-2021}

\RequirePackage{amsthm,amsmath,amsfonts,amssymb}
\RequirePackage[numbers]{natbib}

\startlocaldefs
\endlocaldefs

\usepackage{url}
\usepackage{enumerate}
\usepackage{color}
\usepackage{xr}
\usepackage{amssymb}
\usepackage{amsmath}
\usepackage{latexsym}
\usepackage{epsfig}
\usepackage{textcomp}
\usepackage{amsmath}
\usepackage{graphicx,psfrag,epsf}
\usepackage{enumerate}

\usepackage{fontenc}
\usepackage{url} % not crucial 

\startlocaldefs

 \newcommand{\tr}{\textcolor{red}}

\endlocaldefs

\begin{document}

\begin{frontmatter}
\title{Polyspectral Mean Estimation of General Nonlinear Processes}
%\title{A sample article title with some additional note\thanksref{T1}}
\runtitle{Polyspectral Mean Estimation}
%\thankstext{T1}{A sample of additional note to the title.}

\begin{aug}
\author[A]{\fnms{Dhrubajyoti} \snm{Ghosh}\ead[label=e1,mark]{dg302@duke.edu}},
\author[B]{\fnms{Tucker} \snm{McElroy}\ead[label=e2]{tucker.s.mcelroy@census.gov}},
\and
\author[C]{\fnms{Soumendra} \snm{Lahiri}\ead[label=e3]{s.lahiri@wustl.edu}}
%%%%%%%%%%%%%%%%%%%%%%%%%%%%%%%%%%%%%%%%%%%%%%
%% Addresses                                %%
%%%%%%%%%%%%%%%%%%%%%%%%%%%%%%%%%%%%%%%%%%%%%%
\address[A]{Department of Biostatistics and Bioinformatics, Duke University, \printead{e1}}
\address[B]{Research and Methodology  Directorate, U.S. Census Bureau, \printead{e2}}
\address[C]{Department of Mathematics and Statistics, Washington University in St. Louis, \printead{e3}}
\end{aug}
\begin{abstract}
Higher-order spectra (or polyspectra), defined as the Fourier Transform of a stationary process' autocumulants, are useful in the analysis of nonlinear and non-Gaussian processes. Polyspectral means are weighted averages over Fourier frequencies of the polyspectra, and estimators can be constructed from analogous weighted  averages of the higher-order periodogram (a statistic computed from the data sample's discrete Fourier Transform).  We derive the asymptotic distribution of a class of polyspectral mean estimators,
obtaining an exact expression for the limit distribution that depends on both the given weighting function as well as on higher-order spectra.  Secondly, we
use bispectral means to define a new test of the linear process hypothesis.
 Simulations document  the finite sample properties of the asymptotic results.  Two applications illustrate our results' utility: we test the linear process hypothesis for a Sunspot time series,  and for the Gross Domestic Product
 we conduct a clustering exercise based on bispectral means with different weight functions.
\end{abstract}


\end{frontmatter}

%%%%%%%%%%%%%%%%%%%%%%%%%%%%%%%%%%%%%%%%%%%%%%
%%%% Main text entry area:

\section{Introduction}
\label{sec:intro}
\setcitestyle{authoryear}

Spectral Analysis is an important tool in time series analysis. The spectral density
%, defined as the discrete Fourier transform of the autocovariance function, 
provides important insights into the underlying periodicities of a time series; see 
\citet{brockwell2016introduction} and \citet{mcelroypolitis2020}.
However, when the stochastic process is nonlinear,
%or when the innovations arise from a non-Gaussian process, 
higher-order autocumulants (cf. \citet{brillinger2001time})
can furnish valuable insights that are not 
discernible from the spectral density, as argued in   
%Therein lies the importance of higher order spectra, or polyspectra, as defined by 
\citet{brillinger1967asymptotic}. 
In particular, three- and four-way interactions of a nonlinear process 
can be measured through the third and fourth order polyspectra; 
 see \citet{brillinger1965introduction}, \citet{brillinger1967asymptotic},   \citet{brockett1988bispectral},
 \citet{gabr1988third}, 
 \citet{maravall1983application}, 
 \citet{mendel1991tutorial}, and \citet{berg2009higher}.
Linear functionals of higher-order spectra (here referred to as {\it  polyspectral means}) are often encountered in the
analysis of nonlinear time series, and are therefore a serious object of inference.
Whereas there is an existent literature on estimation of polyspectra
(\citet{brillinger1967asymptotic}, \citet{van1966asymptotic}, \citet{raghuveer1986bispectrum}, \citet{spooner1991estimation}), 
 inference for polyspectral means remains a gap.  This paper aims to fill this gap by developing
 asymptotic distributional results on polyspectral means,  and by developing 
statistical inference methodology for nonlinear time series analysis.

Let $\{X_t\}$ be a $(k+1)^{th}$ order stationary time series, that is,  
 $E [ {| X_t |}^{k+1} ] < \infty$ and \\
 $EX_{t}X_{t+h_1}\ldots X_{t+h_r}
 = EX_{1}X_{1+h_1}\ldots X_{1+h_r}$
 for all integers $t, h_1,\ldots, h_r$,
 and for all $1\leq r \leq k$, 
 for given $k \geq 1$. The most commonly occurring case is $k=1$, corresponding to the 
 class of second order stationary (SOS)  processes. Thus, for a SOS process
 $\{X_t\}$, $EX_t=EX_1$ is  a constant and its cross-covariances depend only on the time-difference 
  $ \cov(X_t, X_{t+h}) = \gamma (h)$ for all
 $t,h\in \bbz$, where 
 $\bbz=\{0,\pm 1,\pm 2, \ldots\}$
 denotes the set of all integers. When the autocovariance 
 function $\ga(h)\equiv \ga_h$ is absolutely summable, $\{X_t\}$ 
 has a spectral density given by $f(\la) = \sum_{h\in \bbz} 
 \ga(h) e^{-\ci h \la}$, $\la\in [-\pi,\pi]$,
 where  $\ci=\sqrt{-1}$. 
 
Important features of a stochastic process  can be extracted with a spectral  mean, which is defined as $\int_{[-\pi, \pi]} f(\lambda)g(\lambda) d\lambda$, where
 %$f(\cdot)$ is the spectral density and 
 $g(\cdot)$ is a weight function. Spectral means can provide  important insights into the time series afforded by different g functions. For example, a spectral mean with  $g(\lambda) = cos(h\lambda)$ corresponds to the lag $h$ autocovariance, whereas $g(\lambda) = \ind(-a,a)$ gives the spectral content in the interval $(-a,a)$, where $\ind(\cdot)$
 denotes the indicator function. Analogous definitions can be made for higher order spectra. 
 The $(k+1)^{th}$ order autocumulant is defined as 
 $\gamma (h_1,\ldots, h_k)
 = \Cum (X_t,X_{t+h_1}, \ldots, X_{t+h_k})$ for all $t$, where
 $\Cum (X,Y,\ldots,Z)$ denotes the cumulant
 of jointly distributed random variables $\{X,Y,\ldots,Z\}$ (cf.  \citet{brillinger2001time})  and where  $h_1, \ldots, h_k$
 are  integer lags.
  If  the corresponding autocumulant function is absolutely summable, then we  define the polyspectra of order $k$ as the Fourier transform of the $(k+1)^{th}$ order autocumulants: 
 \tr{
  $$ 
  f_k(\underline{\lambda}) = \sum_{\underline{h} \in \bbz^k} \gamma (\underline{h}) e^{-\iota \underline{h}'\underline{\lambda}} \, ,
  $$
  }
  where  $\underline{h} = {[ h_1, \ldots, h_k]}^{\prime}$,
  $\underline{\lambda} = {[\lambda_1, \ldots, \lambda_k ]}^{\prime}$
  and $A'$ denotes the transpose of a matrix $A$. 
  %In the case that $k=1$,the autocumulants are autocovariances, and the polyspectrum is the spectral density.
A polyspectral mean
with weight function $g: [-\pi,\pi]^k \raw \mathbb{C}$
is defined as 

\beq
M_g(f_k)=\int_{[-\pi,\pi]^k} f_k(\underline{\lambda})g(\underline{\lambda})d\underline{\lambda}.
\label{psp-m}
\eeq
Thus, spectral means correspond to the case $k=1$, in which case we drop the subscript $k$ (and write $f_1=f$). 
 Polyspectral means give us important information about a time series that can  not be obtained from the spectral distribution, e.g.,  when the time series is nonlinear or the innovation process is non-Gaussian. 
%We have observed in [Our First Paper] that it is possible to gain prediction accuracy by using quadratic predictor rather than linear predictor when there is nonlinearity in a time series. 
%Hence it is of significant importance to identify whether there is nonlinearity in a time series in order to use accurate prediction procedures.  
%
Additionally, we can extract several features from a time series through the polyspectral mean corresponding to different weight functions. These features can play an important role in identifying (dis-)similarities among 
different  time series
%and as a result, polyspectral means have been widely used by 
%researchers in applied sciences and engineering 
(cf. \citet{chua2010application}, \citet{picinbono1999polyspectra}, \citet{dubnov1996polyspectral}, \citet{vrabie2003spectral}). 

\begin{figure}[htbp]
    \centering
    \includegraphics[width = 0.6\textwidth]{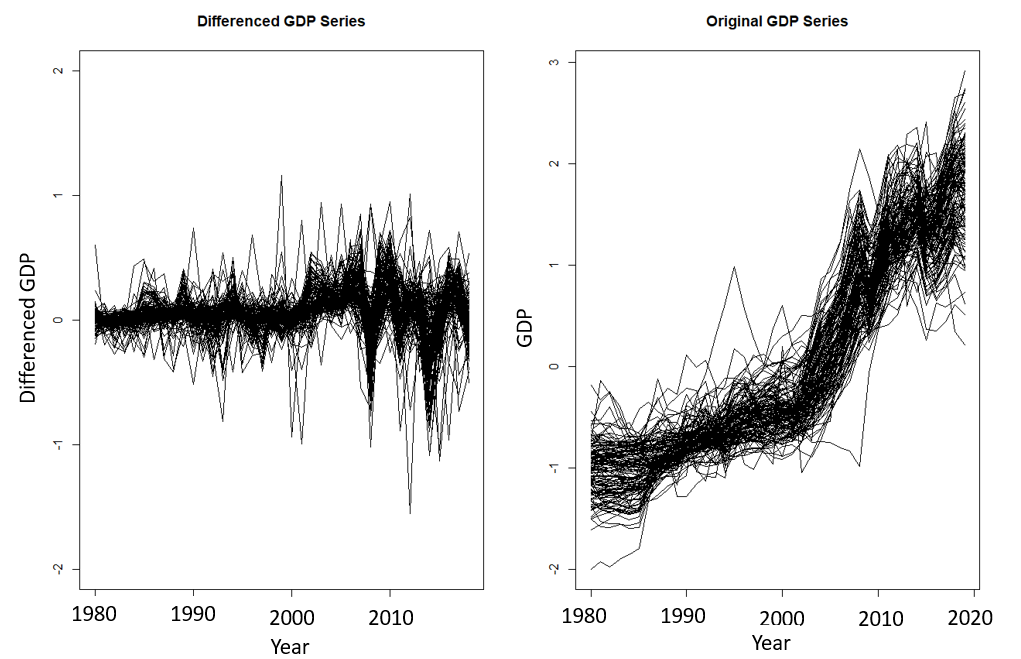}
    \caption{Differenced GDP (left panel) and raw GDP (right panel)  for 136 countries, annual
    1980-2020. There is a common nonlinear trend
    (right panel). 
    %There are three main types of trends in the data, continuously increasing, increasing with a downward trend toward the end, and increasing with a downward trend followed by an increasing trend at the end.
    }
    \label{fig:gdp_figure}
\end{figure}

For an illustrative example, 
we consider the annual Gross Domestic Product (GDP)  data from 136 countries over 40 years. Our goal is to 
obtain clustering of the countries based on patterns of their GDP growth rates. 
% economic activities as captured by the GDP data.
The right panel of Figure \ref{fig:gdp_figure}   gives the raw time series and the left panel gives the differenced and scaled versions of the GDP of different countries. As a part of exploratory data 
analysis, we tested for Gaussianity and   linearity;
 a substantial proportion of the time series are non-Gaussian, and some of them are nonlinear. Hence, using  spectral means alone to capture relevant
 time series features 
may  be inadequate. Here we use (estimated) higher order polyspectral means with different
weight functions to elicit salient features  of the GDP data, in order 
to capture possible nonlinear  dynamics of  the economies. See Section \ref{sec:data}
for more details.

 Estimation  of spectral and polyspectral 
means can be carried out 
using the periodogram and its higher order versions, which are defined in terms of the discrete Fourier transform of
a length $T$ sample $\{X_1,\ldots, X_T\}$ 
 (described in Section \ref{sec:expl}). However, distributional 
properties of the polyspectral mean 
estimators are not
very well-studied. 
%In an important work, 
 \citet{dahlhaus1985asymptotic} proved the asymptotic normality of spectral means, i.e., the case $k=1$; however, in the general case ($k > 1$),
 there does not seem to be
 any work on the asymptotic distribution of polyspectral mean estimators. One of the major contributions of this 
 paper is to establish asymptotic normality results for the polyspectral means of a general order. We develop some  nontrivial combinatorial arguments involving higher-order autocumulants to show that, 
 under mild conditions, the estimators of the  $k^{th} $ order polyspectral mean 
 parameter $M_g(f_k)$ 
 in \eqref{psp-m} are asymptotically normal.
 We also obtain an explicit expression for the asymptotic  variance, which is shown to 
 depend on certain 
 polyspectral means of order $2k+1$.
 This result agrees with the 
 known results on  spectral means when 
 specialized to the case $k=1$, where 
 the limiting variance 
 is known to involve
 the trispectrum $f_{2k+1}=f_3$.
 In particular, the results of this paper 
 provide a unified way  to construct 
 studentized versions of spectral
 mean estimators and  higher order polyspectral mean estimators, which  can be
 used to carry out 
 large sample statistical inference on 
 polyspectral mean parameters of any 
{\it arbitrary}
 order $k\geq 1$.

 The second major  contribution of the paper is the development of a
 new test procedure to assess the linearity assumption on a time series. 
 Earlier work on the problem is typically based on the squared modulus of
 the estimated bispectrum (\citet{rao1980test}, \citet{chan1986tests}, \citet{berg2010bootstrap}, \citet{barnett2004time}). In contrast, here 
 we make a key observation that under the linearity hypothesis (i.e., the 
 null hypothesis) the ratio of the bispectrum to a suitable 
 power transformation of the spectral density must be a constant
 at {\it all} frequency pairs in $[-\pi,\pi]^2$. This observation
 allows us to construct a set of higher order polyspectral means that must be 
 zero when the process is linear. On the other hand, for a nonlinear 
 process, not all of these spectral means can be equal to zero. 
 Here we exploit this fact and 
 develop a new test statistic that can be used to
 test the significance of these 
 polyspectral means.
 We also  derive the asymptotic distribution of the 
test statistic  under the null hypothesis, and provide
the critical values needed for
calibrating the proposed test.

The rest of the paper is organized as follows. Section
\ref{sec:expl} covers some background material, and provides some examples 
of polyspectra and polyspectral means.
Section 
\ref{sec:results} gives the regularity conditions and the 
asymptotic properties of the polyspectral means of
a general order. The linearity test is described in Section 
\ref{sec:application}.
Simulations are presented in Section \ref{sec:sim}, and
two data applications, one involving the 
Sunspot data and the 
other involving the GDP growth rate 
(with a discussion of clustering time series
via bispectral means), are presented in Section \ref{sec:data}.
The proofs of the theoretical results are given in Section \ref{sec:proof}.
% we examine the Gross Domestic Product (GDP) growth rate  to cluster countries around the world.
% based on the pattern of progression of the GDP over past 40 years.

\section{Background and examples}
\label{sec:expl}
\setcounter{equation}{0}
% \vspace{-0.8cm}

First, we will define estimators of the polyspectral mean 
 $M_g\left(f_k \right)
 %= \int_{[-\pi, \pi]^k} f_k\left(  \underline{\lambda} \right)g\left( \underline{\lambda} \right) d\underline{\lambda}
 $ defined  in \eqref{psp-m}, 
 corresponding to a given weight function $g$. 
 Following \citet{brillinger1967asymptotic}, we define an estimator of 
 $M_g\left(f_k \right)$
 by replacing $f_k$ by the $k^{th}$ order periodogram.
Specifically, given a sample $\{ X_1, \ldots, X_T \}$ from the  mean zero time series $\{ X_t \}$,
let  $d(\lambda) = \sum_{t=1}^T X_t e^{-\iota \lambda t}$ denote  the Discrete Fourier Transform (DFT)
 of the sample, where  $\lambda \in (-\pi, \pi]$
 %%%%%%%%%% C: SNL
 % We need to restrict the number of Fourier
 %freq to $T$.
 % This requires 1 less integer in the 
 % allowable set of $\ell$  below in the even case.
 % Restricting $\la\in (-\pi, \pi]$ 
 % in place of  $\la\in [-\pi, \pi]$  fixes this %problem
 %%%%%%%%%%%%%%%
 and  $\iota = \sqrt{-1}$.
 Such a frequency may be a Fourier frequency,
 taking the form $2 \pi \ell/T$ for some
 integer $\ell \in \{ -\floor{\frac{T}{2}},
  \ldots, \floor{\frac{T}{2}} \}$, in which
  case we write $\tilde{\lambda}$.  Here and in the following,
  $\floor{x}$ denotes the integer part of a real number $x$. 
  A $k$-vector of frequencies is denoted
  by $\underline{\lambda} = {[\lambda_1, \ldots,  \lambda_k]}'$; if these are all Fourier  frequencies, then we write $\tilde{\underline{\lambda}}$.
 Then the $k^{th}$ order periodogram is defined as
 \[
 \hat{f}_k(\underline{\lambda})= T^{-1}d (\lambda_1 )\cdots d(\lambda_k )d 
 (- [\underline{\lambda}] ),
 \] 
 where  $[\underline{\lambda}] $ is a shorthand for $\sum_{j=1}^k \lambda_j$.
Next, let $\sum_{\tilde{\underline{\lambda}}}$ denote
 a shorthand for a summation over
 vector $\underline{\lambda}$ that consist
 of Fourier frequencies.
% $-\floor{\frac{T}{2}} \leq j_1,\ldots,j_k \leq \floor{\frac{T}{2}} $.  
Then we  define the estimator of the  polyspectral mean  $M_g\left(f_k \right)$ as
\begin{align}
    \nonumber \widehat{M_g(f_k)} &\equiv  (2 \pi)^k T^{-k}\sum_{\tilde{\underline{\lambda}}} \hat{f}_k(\tilde{\underline{\lambda}})
    g\left(\tilde{\underline{\lambda}} \right)\Phi(\tilde{\underline{\lambda}}) \\
    &= (2 \pi)^k T^{-k-1} \sum_{\tilde{\underline{\lambda}}} d\left(\tilde{\lambda_1} \right)\cdots d\left(\tilde{\lambda_k} \right)d\left(-  [\tilde{\underline{\lambda}} ] \right)g\left(\tilde{\underline{\lambda}}\right)
    \Phi(\tilde{\underline{\lambda}}), 
    \label{eq:polyHat}
\end{align}
where $\Phi(\underline{\lambda})$ is an indicator function that is zero  
if and only if  $\underline{\lambda}$  lies in any proper linear sub-manifold of the $k$-torus ${[-\pi, \pi]}^k$. That is, $\Phi(\underline{\lambda})$
 is one if and only if for all  subset $\lambda_{i_1}, \ldots, \lambda_{i_m}$ of $\underline{\lambda}$ ($1\leq m\leq k$),
 $\sum_{j=1}^m \lambda_{i_j} \not\equiv 0 \textrm{ (mod } 2 \pi)$. 
 We further assume that the weighting function $g$ is continuous in $\underline{\lambda}$
 (except possibly on a set of Lebesgue measure
 zero; cf. Section 2.4, \citet{AthreyaLahiri2006} 
 ), and satisfies the symmetry condition
\begin{equation}
    \label{eq:sym-g}
      g(\underline{\lambda}) = \overline{g(-\underline{\lambda})}
\end{equation}  
 for all $\underline{\lambda}$.  This
 condition holds for all polyspectra, 
 and holds more broadly for any function that
 has a multi-variable Fourier expansion
 with real coefficients, thereby
 encompassing most functions of interest
 in applications.  Moreover, 
 this assumption 
 ensures that the corresponding
 polyspectral mean is real-valued, since
\[
 \overline{M_g (f_k)}  = \int_{[-\pi,\pi]^k} \overline{f_k(\underline{\lambda}) } 
 \overline{ g(\underline{\lambda}) }
  \, d\underline{\lambda} =
  \int_{[-\pi,\pi]^k} 
  f_k(-\underline{\lambda}) 
  g(-\underline{\lambda})   \, d\underline{\lambda} = M_g (f_k),
\]
using  (\ref{eq:sym-g}) and
a change of variable $-\underline{\lambda} 
\mapsto \underline{\lambda}$. 
We will also  assume that $\int_{\underline{\lambda}} \lvert g\left( \underline{\lambda} \right) \rvert < \infty$,
 which screens out perverse cases of
 that are of little practical interest.

 Note that the value of the polyspectral mean (\ref{psp-m}) is the same whether or not we
 restrict integration to sub-manifolds of
 the torus,   since such sets have Lebesgue measure zero.
 %; hence any estimators must exclude $\underline{\lambda}$ that are contained in such sub-manifolds. 
  We need to introduce the indicator function, since \citet{brillinger1967asymptotic} established that the $k^{th}$ order periodogram diverges if the Fourier frequencies lie in 
  such sub-manifolds.  This is a significant issue. Indeed,
  the estimator of the polyspectral density 
   in \citet{brillinger1965introduction} 
  is also given by a kernel-weighted average of the $k^{th}$ order periodogram, where the average is taken by avoiding the sub-manifolds. 
  The polyspectral density estimator
involves a bandwidth term, which in general makes the convergence rate much slower. Since we are working with polyspectral means, it is not necessary to smooth the $k^{th}$ order periodogram, and hence 
we can ignore the bandwidth problem and 
focus directly  on regions that avoid the 
 sub-manifolds. 
 
\begin{figure}[htbp]
    \centering
    \includegraphics[width = 0.3 \textwidth]{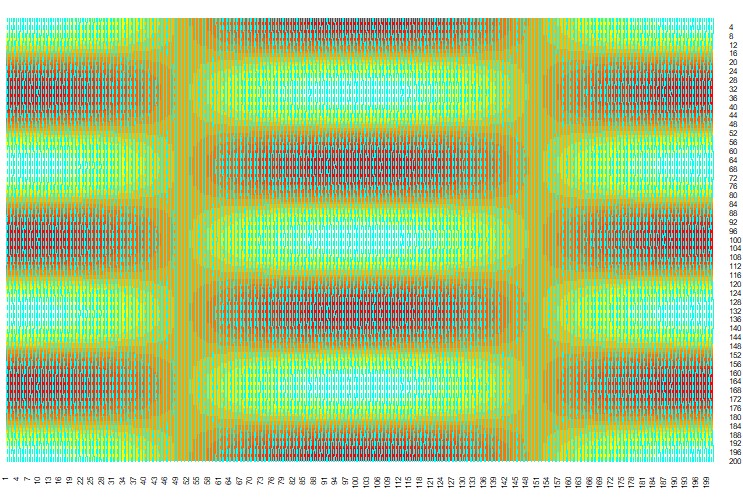}
    \includegraphics[width = 0.3 \textwidth]{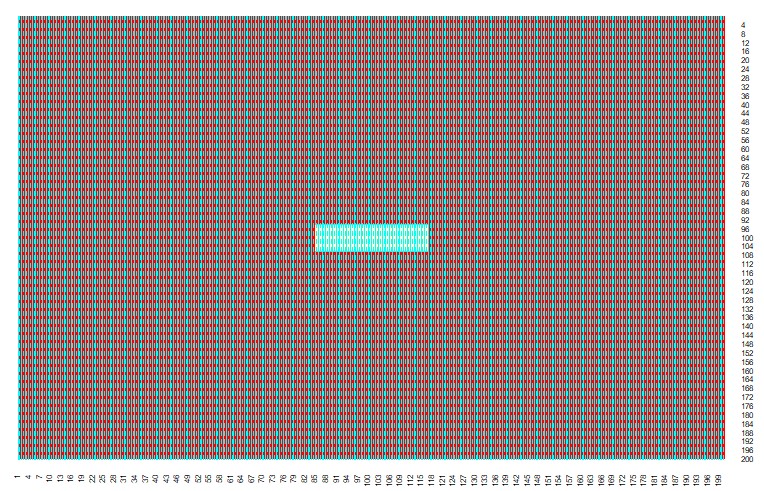}
    \includegraphics[width = 0.3 \textwidth]{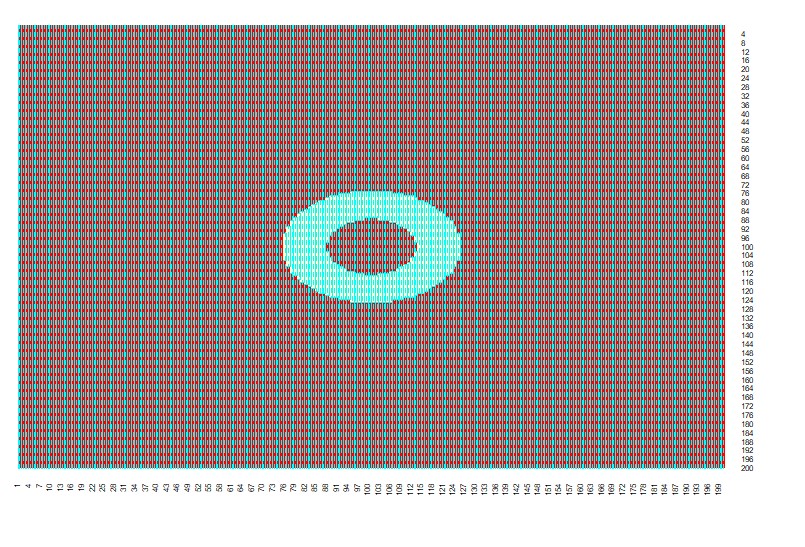}
    \caption{Heatmap for different g functions: (left panel) $g(\underline{\lambda}) = cos(2 \lambda)cos(3 \lambda)$, (center panel) $g(\underline{\lambda}) = \mathbb{I}(\lambda_1 \leq 0.2)\mathbb{I}(\lambda_2 \leq 0.2)$, (right panel) $g(\underline{\lambda}) = \mathbb{I}( 0.1 \leq \lambda_1^2 + \lambda_2 ^2 \leq 0.2)$.}
    \label{fig:heatmap_weight_function}
\end{figure}
  
  Some illustrations of
 weighting functions  are given in Figure \ref{fig:heatmap_weight_function}. We
 also  provide some specific examples of polyspectral means below,  demonstrating how diverse features of
 a nonlinear process can be extracted with different weighting functions.  (Recall that we assumed that the process has mean zero, which assumption is used throughout.)
%%%
\begin{example}[Example 1:]
\label{exam:g1}
    Our first example, which was mentioned
    in the introduction, is already familiar:
    let $g(\lambda_1) = {(2 \pi)}^{-1}
    \exp \{ \iota \lambda_1 h_1 \}$, so that
$ 
M_g (f_1) = \int_{-\pi}^{\pi} f_1 (\lambda_1)
 g(\lambda_1) d\lambda_1 = \gamma (h_1)$,
where $\gamma (\cdot) $ is the second
order autocumulant function, i.e., 
 the autocovariance function.  Its classical
 estimator  is the sample autocovariance
$ \hat{\gamma} (h_1) = T^{-1} \sum_t X_t X_{t+h_1}$,
 where the summation is over all $t$ such that
$t, t+h_1 \in \{ 1, 2, \ldots, T \}$ (the sampling times).
This agrees with our proposed polyspectral mean
estimator, viz.
$ 
    \widehat{M_g (f_1)}
     = T^{-2} \sum_{\tilde{\lambda_1}}
     d (\tilde{\lambda_1}) d (-\tilde{\lambda_1})
     g( \tilde{\lambda_1}) = \hat{\gamma} (h_1) $,
using the fact that
\begin{equation}
    \label{eq:ft-elem}
    \sum_{\tilde{\lambda_1}} \exp \{ \iota \tilde{\lambda_1} s \} = T \, \ind_{ \{ s= 0 \} }.
\end{equation}
Note that in our
expression for $\widehat{M_g (f_1)}$ there
are no sub-manifolds of the 1-torus.
\end{example}

% \begin{example}[Example 2:]
% \label{exam:g0}
    Example 1
    % \ref{exam:g1} 
    can be extended to third-order autocumulants via the weight function $g(\lambda_1, \lambda_2) = {(2 \pi)}^{-2} \exp \{ \iota \lambda_1 h_1 + \iota \lambda_2 h_2 \}$.  Then it follows that
\begin{equation}
    M_g(f_2) = \int_{-\pi}^\pi \int_{-\pi}^\pi f_2(\underline{\lambda})g(\underline{\lambda})d \underline{\lambda} = \gamma(h_1, h_2), 
\end{equation}
        % M_g(f_2) = \int_{-\pi}^\pi \int_{-\pi}^\pi f_2(\underline{\lambda})g(\underline{\lambda})d \underline{\lambda} 
        % &= (4 \pi^2)^{-1}\int_{-\pi}^\pi \int_{-\pi}^\pi \sum_{k_1, k_2 \in \mathbb{Z}} \gamma(k_1, k_2)e^{\iota \lambda_1(h_1 - k_1) + \iota \lambda_2 (h_2 - k_2)} d \lambda_1 d\lambda_2 \\
        % = \gamma(h_1, h_2), 
 % \]
        where $\gamma(\cdot, \cdot)$ is the third order autocumulant function.  The
        classical estimator
        of this autocumulant function is given by
    $    \hat{\gamma}(h_1, h_2) = T^{-1}\sum_{t} X_{t}X_{t + h_1} X_{t + h_2},    $
    where the summation is over all $t$ such that
    $t, t+h_1, t+h_2 \in \{ 1, 2, \ldots, T\}$.
 We next show that this is asymptotically
  equivalent (i.e., the same up to terms tending
  to zero in probability) to 
  $\widehat{M_g(f_2)}$.  First, note that
  the sub-manifolds of the 2-torus that we
   must exclude are simply the sets where
   either $\tilde{\lambda}_1$ or $\tilde{\lambda}_2$
   equals zero; 
   %%%%%%%%%%%%%%%%%%%%%%%%%%%%%%%%%%%%%%%%%%
   %%SNL
  % {\bf C: How about excluding 
  % $\tilde{\lambda}_1+\tilde{\lambda}_2=0$ ? }
  %%%%%%%%%%%
  %% The final contribution is $O(T^{(k-1)})$
  %% making the contribution negligible
  %%%%%%%%%%%%%%%%%%%%%%%%%%%%%%%%%%%%%%%%%%
   hence
%     \[
%       \widehat{M_g(f_2)} =  T^{-3} \sum_{\tilde{\lambda_1}, \tilde{\lambda_2} \neq 0} d(\tilde{\lambda_1})d(\tilde{\lambda_2}) d(-\tilde{\lambda_1} - \tilde{\lambda_2})
%       e^{\iota \tilde{\lambda_1} h_1 + \iota \tilde{\lambda_2} h_2}
% \]      
\begin{align*}
  \widehat{M_g(f_2)} &=  T^{-3} \sum_{\tilde{\lambda_1}, \tilde{\lambda_2} \neq 0} d(\tilde{\lambda_1})d(\tilde{\lambda_2}) d(-\tilde{\lambda_1} - \tilde{\lambda_2})
      e^{\iota \tilde{\lambda_1} h_1 + \iota \tilde{\lambda_2} h_2} \\
    &= T^{-3}\sum_{\tilde{\lambda_1}, \tilde{\lambda_2}}\sum_{t_1, t_2, t_3}X_{t_1}X_{t_2}X_{t_3}e^{-\iota \tilde{\lambda_1} t_1}e^{-\iota \tilde{\lambda_2} t_2}e^{\iota (\tilde{\lambda_1} + \tilde{\lambda_2})t_3}e^{\iota \tilde{\lambda_1} h_1}e^{\iota \tilde{\lambda_2} h_2} \\
        &- T^{-3}\sum_{\tilde{\lambda_1}}\sum_{t_1, t_2, t_3}X_{t_1}X_{t_2}X_{t_3}e^{\iota \tilde{\lambda_1} (t_3 - t_1 + h_1)} - T^{-3}\sum_{\tilde{\lambda_2}}\sum_{t_1, t_2, t_3}X_{t_1}X_{t_2}X_{t_3}e^{\iota \tilde{\lambda_2} (t_3 - t_2 + h_2)} \\
        & + T^{-3} \sum_{t_1, t_2, t_3} 
         X_{t_1} X_{t_2} X_{t_3}.
\end{align*}
 Note that the last term is simply the cube
 of the sample mean $\overline{X}$,
 which tends to zero
 in probability.  For the other terms,
 we can apply (\ref{eq:ft-elem}) to obtain
%
 %%%\begin{align*}  
  $  \widehat{M_g(f_2)}  = 
    T^{-1}\sum_{t_1}X_{t_1}X_{t_1 + h_1}X_{t_1 + h_2}  
   - T^{-1}\sum_{t_2}X_{t_2} T^{-1}\sum_{t_1}X_{t_1}X_{t_1 - h_1} - T^{-1}\sum_{t_1}X_{t_1} T^{-1}\sum_{t_2}X_{t_2}X_{t_2 - h_2}   
   +    { ( \overline{X})}^3
$,
   %%%%%%%%%%%%%%%%%%%%%%%%%%%%%%%%%%%%%%%
  %      &\sim T^{-3} \sum_{\tilde{\lambda_1}, \tilde{\lambda_2}}\sum_{t_1, t_2, t_3}X_{t_1}X_{t_2}X_{t_3}e^{-\iota \tilde{\lambda_1}(t_1 - t_3 - h_1) - \iota \tilde{\lambda_2} (t_2 - t_3 - h_2)} \\
%      & + (2 \pi)^{-2}T^{-1}\sum_{t_1,t_2, t_3} X_{t_1}X_{t_2}X_{t_3} \int_{-\pi}^\pi \int_{-\pi}^\pi e^{\iota\tilde{\lambda_1}(t_3 - t_1 + h_1) + \iota \tilde{\lambda_2} (t_3 - t_2 + h_2)}d\lambda_1 d\lambda_2 \\
 %    & \quad \quad - T^{-2}(2\pi)^{-1}\sum_{t_1, t_2, t_3}X_{t_1}X_{t_2}X_{t_3}\int_{-\pi}^{\pi}e^{\iota \tilde{\lambda_1} (t_3 - t_1 + h_1)}d\lambda_1  \\ \displaybreak
%     & \quad \quad  - T^{-2}(2\pi)^{-1} \sum_{t_1, t_2, t_3}X_{t_1}X_{t_2}X_{t_3}\int_{-\pi}^{\pi}e^{\iota \tilde{\lambda_2} (t_3 - t_2 + h_2)}d\lambda_2 \\ 
%%%% \end{align*}
 which is asymptotically equivalent to
  $\hat{\gamma}(h_1, h_2)$.

\begin{example}[Example 2:]
\label{exam:g2}
For a zero mean process, the autocumulants and automoments of order 2 and 3 coincide, and therefore
Example 1 shows that the spectral mean estimator 
is a consistent estimator of these autocumulants/automoments.
A similar (but more involved) calculation
is required to show that a similar result holds for 
higher order autocumulants.
As an example, we now consider the fourth order
autocumulant function (which is different  from  
the fourth order automoment function). 
To that end,  we consider   the trispectral mean with the weight function 
    $ g(\lambda_1, \lambda_2, \lambda_3) = {(2 \pi)}^{-3}\exp \{ \iota (\lambda_1 h_1 + \lambda_2 h_2 + \lambda_3 h_3) \}$.  Then, letting $\eta (k_1, k_2, k_3)$ denote the order 3 auto-moment,
    and using the identity
$ \gamma(k_1, k_2,k_3)= \eta(k_1, k_2,k_3) - \gamma(k_1)\gamma(k_3 - k_2) -\gamma(k_2)\gamma(k_3 - k_1) - \gamma(k_3)\gamma(k_2 - k_1)$, we get
  \begin{align*}
  % \label{eq:}
        M_g(f_3) &
        %= \int_{-\pi}^\pi \int_{-\pi}^\pi \int_{-\pi}^\pi f_3(\underline{\lambda})g(\underline{\lambda})d \underline{\lambda} \\
        %&
        = (8 \pi^3)^{-1}\int_{-\pi}^\pi \int_{-\pi}^\pi \int_{-\pi}^\pi \sum_{k_1, k_2,k_3 \in \mathbb{Z}} \gamma(k_1, k_2,k_3) \\
        & \quad \quad \quad e^{\iota \lambda_1(h_1 - k_1) + \iota \lambda_2 (h_2 - k_2) + \lambda_3 (h_3 - k_3)} d \lambda_1 d\lambda_2 d\lambda_3 \\
     %   &= (8 \pi^3)^{-1}\int_{-\pi}^\pi \int_{-\pi}^\pi \int_{-\pi}^\pi \sum_{k_1, k_2,k_3 \in \mathbb{Z}} \Big\{\eta(k_1, k_2,k_3) - \gamma(k_1)\gamma(k_3 - k_2)  \\
     %   &\quad \quad \quad -\gamma(k_2)\gamma(k_3 - k_1) - \gamma(k_3)\gamma(k_2 - k_1) \Big\} \\
     %   & \quad \quad \quad e^{\iota \lambda_1(h_1 - k_1) + \iota \lambda_2 (h_2 - k_2) + \lambda_3 (h_3 - k_3)} d \lambda_1 d\lambda_2 d\lambda_3 \\
        &= \eta(h_1, h_2, h_3) - \gamma(h_1)\gamma(h_3 - h_2) - \gamma(h_2)\gamma(h_1 - h_3) - \gamma(h_3)\gamma(h_2 - h_1).
        \end{align*}
  %      where $\gamma(\cdot, \cdot)$ is the third order automoment
        %autocumulant 
   %     function.
Setting $S = \{ \lambda_1, \lambda_2, \lambda_3\neq 0 : \lambda_1 + \lambda_2 \neq 0, \lambda_2 + \lambda_3 \neq 0, \lambda_1 + \lambda_3 \neq 0 \}$, we obtain
\begin{align*}
        \widehat{M_g(f_3)} &= (2 \pi)^{3} T^{-4} \sum_{\underline{\tilde{\lambda}} \in S} d(\tilde{\lambda_1})d(\tilde{\lambda_2})d(\tilde{\lambda_3})d(-\tilde{\lambda_1} - \tilde{\lambda_2} -\tilde{\lambda_3})g(\tilde{\lambda_1}, \tilde{\lambda_2}, \tilde{\lambda_3}) \\
        &= T^{-4} \sum_{\underline{\tilde{\lambda}} \in S}   \sum_{t_1, t_2, t_3, t_4} X_{t_1}X_{t_2}X_{t_3}X_{t_4}e^{-\iota \tilde{\lambda_1} (t_1 - t_4 - h_1) - \iota \tilde{\lambda_2} (t_2 - t_4 - h_2) -\iota \tilde{\lambda_3} (t_3 -t_4 - h_3)}. 
        %  & \quad \quad \quad \qquad \qquad \qquad  cos(\lambda_1h_1)cos(\lambda_2h_2)cos(\lambda_3h_3) d\lambda_1 d\lambda_2  \\
        % &- (2 \pi)^{-2} T^{-1} T^{-2} \sum_{\lambda_1}\sum_{t_1. t_2, t_3}X_{t_1}X_{t_2}X_{t_3}e^{-\iota \lambda_1 (t_1-t_3) }  cos(\lambda_1 h_1)d\lambda_1 \\
        % &- (2 \pi)^{-2} T^{-1} T^{-2} \sum_{\lambda_2}\sum_{t_1. t_2, t_3}X_{t_1}X_{t_2}X_{t_3}e^{-\iota \lambda_2 (t_2-t_3) }  cos(\lambda_2 h_2)d\lambda_2 \\
        % &=   T^{-1}\sum_{t_1, t_2, t_3} X_{t_1}X_{t_2}X_{t_3} \int_{-\pi}^\pi  \int_{-\pi}^\pi cos(\lambda_1 h_1)e^{\iota \lambda_1(t_3 - t_1)} cos(\lambda_2 h_2)e^{\iota \lambda_2(t_3 - t_2)} d\lambda_1 d\lambda_2 \\
        % &- (2\pi)^{-1}(\Bar{X_t}) T^{-1}\sum_{t_1}X_{t_1}X_{t_1+h_1} - (2\pi)^{-1}(\Bar{X_t}) T^{-1}\sum_{t_1}X_{t_1}X_{t_1+h_2} \\
        % &= T^{-1} \sum_{t_1} X_{t_1}X_{t_1 + h_1 - h_2}X_{t_1 + h_1} - (2\pi)^{-1}(\Bar{X_t}) (\hat{\gamma}(h_1) + \hat{\gamma}(h_2))\\
        % &=  \hat{\gamma}(h_1 - h_2, h_1)  + o(1)
    \end{align*}
    If $\lambda_i=0$ for some $i$, then we will obtain a $\bar{X}$ term, which will be asymptotically zero under the zero mean assumption. Consider the case for $\lambda_1 + \lambda_2 = 0$:
   \quad $ 
   T^{-4} \sum_{\lambda_1, \lambda_2} \sum_{t_1, t_2, t_3, t_4} X_{t_1}X_{t_2}X_{t_3}X_{t_4} \times$ 
   $e^{-\iota \lambda_1(t_1-t_2-h_1+h_2)}e^{-\iota \lambda_2 (t_3-t_4-h_3)}  
        = \hat{\gamma}(h_1-h_2)\hat{\gamma}(h_3)+o_P(1)
    $.  
    Hence, the estimate takes the form
   \[
         \widehat{M_g(f_3)} = \hat{\eta}(h_1,h_2,h_3) - \hat{\gamma}(h_1)\hat{\gamma}(h_2-h_3) -\hat{\gamma}(h_2)\hat{\gamma}(h_3 - h_1) - \hat{\gamma}(h_3)\hat{\gamma}(h_2-h_1) +o_P(1).
  \]
Each of the terms above is a consistent estimator and therefore,  $\widehat{M_g(f_3)}$ is consistent for $M_g(f_3)$.    
\end{example}

\begin{example}[Example 3:]
\label{exam:g5}
The weighted center of Bispectrum (WCOB)
\citep{zhang2000bispectrum} is a widely used 
feature in nonlinear time series analysis
(cf. 
\cite{chua2010application, acharya2012application, yuvaraj2018novel, martis2013application})
which is defined using ratios of 
polyspectral means. Specifically, we can define the WCOB as $(c_1, c_2) $ 
    (in the $(\lambda_1, \lambda_2)$-torus), where 

    % \tr{
\begin{align*}
 c_1 &= \frac{\int_{-\pi}^{\pi} \int_{-\pi}^{\pi} \lambda_1 f_2(\lambda_1, \lambda_2)d\lambda_1d\lambda_2}{\int_{-\pi}^{\pi}  \int_{-\pi}^{\pi} f_2(\lambda_1, \lambda_2)d\lambda_1d\lambda_2} ~\mbox{and}~
        c_2 = \frac{\int_{-\pi}^{\pi} \int_{-\pi}^{\pi} \lambda_2 f_2(\lambda_1, \lambda_2)d\lambda_1d\lambda_2}{ \int_{-\pi}^{\pi}  \int_{-\pi}^{\pi} 
 f_2(\lambda_1, \lambda_2)d\lambda_1d\lambda_2}
\end{align*}
% }

A discrete analog of WCOB was proposed in \citep{zhang2000bispectrum} which was defined as follows:
    \begin{align*}
        c_1 &= \frac{\sum_{\lambda_1, \lambda_2} \lambda_1 f_2(\lambda_1, \lambda_2)}{\sum_{\lambda_1, \lambda_2} f_2(\lambda_1, \lambda_2)} ~\mbox{and}~
        c_2 = \frac{\sum_{\lambda_1, \lambda_2} \lambda_2 f_2(\lambda_1, \lambda_2)}{\sum_{\lambda_1, \lambda_2} f_2(\lambda_1, \lambda_2)}.
    \end{align*}
%    where $\lambda_1$ and $\lambda_2$ are frequencies in the desired region.

% \tr{A continuous analog of the above example can be given as:
% \begin{align*}
%  c_1 &= \frac{\int_{-\pi}^{\pi} \int_{-\pi}^{\pi} \lambda_1 f_2(\lambda_1, \lambda_2)d\lambda_1d\lambda_2}{\int_{-\pi}^{\pi}  \int_{-\pi}^{\pi} f_2(\lambda_1, \lambda_2)d\lambda_1d\lambda_2} ~\mbox{and}~
%         c_2 = \frac{\int_{-\pi}^{\pi} \int_{-\pi}^{\pi} \lambda_2 f_2(\lambda_1, \lambda_2)d\lambda_1d\lambda_2}{ \int_{-\pi}^{\pi}  \int_{-\pi}^{\pi} 
%  f_2(\lambda_1, \lambda_2)d\lambda_1d\lambda_2}
% \end{align*}}

The quantities $c_1$ and $c_2$ are
ratios of polyspectral means, 
 where the numerators have  weight functions  $g(\lambda_1, \lambda_2) = \lambda_1$
 and $g(\lambda_1, \lambda_2) = \lambda_2$,
 respectively for the cases of 
 $c_1$ and $c_2$; each denominator has a constant weight function. The estimate of the bispectral mean in the numerator of $c_1$ is
\[
 \widehat{M_g(f)} = (2 \pi)^2 T^{-3} \sum_{\tilde{\lambda_1}, \tilde{\lambda_2} \neq 0}  \tilde{\lambda_1} d(\tilde{\lambda_1})d(\tilde{\lambda_2})d(\tilde{-\lambda_1} - \tilde{\lambda_2}).
\]
Estimators of the other polyspectral means 
in $(c_1, c_2)$ can be defined similarly.

    % \begin{align*}
    %     \widehat{M_g(f)} &= (2 \pi)^2 T^{-2} \sum_{\lambda_1, \lambda_2 \neq 0}  T^{-1} \lambda_1 d(\lambda_1)d(\lambda_2)d(-\lambda_1) \\
    %     & \sim T^{-1} \int_{-\pi}^\pi \int_{-\pi}^\pi \sum_{t_1, t_2, t_3} X_{t_1} X_{t_2} X_{t_3} \lambda_1 e^{-\iota \lambda_1 t_1 - \iota \lambda_2 t_2 + \iota (\lambda_1 + \lambda_2)t_3} \\
    %     &= \iota T^{-1} \sum_{t_1, t_2, t_3} X_{t_1}X_{t_2}X_{t_3} \left\{\frac{2 sin (\pi(t_3 - t_1)) - \pi(t_3 - t_1)cos(\pi(t_3 - t_1))}{(t_3 - t_1)^2} \frac{2 sin \pi(t_3 - t_1)}{(t_3 - t_1)}  \right\}
    % \end{align*}
%
   % As we can see, this is a complex number, and hence it is much more prudent\footnote{DJ: calculation is wrong, and should depend on the region anyway.  I commented it out until it
    %can be fixed. Integral over $\lambda_2$ is
    %zero unless $t_2 = t_3$.  Why is it more
    %prudent, what's wrong with a complex
    %polymean?} to use the weight function %$g(\lambda_1, \lambda_2) = \lvert \lambda_1 \rvert$, which is also a symmetric function and hence falls under the assumption of our theory. If we use this as a weight function, the bispectral mean estimate becomes:
    %\begin{align*}
  %      T^{-1} \sum_{t_1, t_2, t_3} X_{t_1}X_{t_2}X_{t_3} \left\{ \frac{(2 (\pi (t_3 - t_1) sin(\pi (t_3 - t_1)) + cos(\pi (t_3 - t_1)) - 1))}{(t_3 - t_1)^2} \frac{2 sin \pi(t_3 - t_1)}{(t_3 - t_1)}  \right\}
%    \end{align*}
\end{example}

\begin{example}[Example 4:]
\label{exam:g3}
 Another type of polyspectral mean is
 obtained by examining the total spectral
 content in a band.  In the $k=1$ case, 
 this has been used to examine the contribution
 to process' variance due to a range
 of stochastic dynamics, such as business
 cycle frequencies.  Analogously for $k=2$,
 we define   $g(\lambda_1, \lambda_2) = \ind(\lambda_1 \in (- \mu_1, \mu_1), \lambda_2 \in (-\mu_2, \mu_2))$ for $\mu_1, \mu_2 > 0$; the corresponding $M_g(f_2)$ measures the bispectral content in the rectangular region $(-\mu_1, \mu_1) \times (-\mu_2, \mu_2)$.
    %  \begin{align*}
    %      M_g(f_2) &= \int_{-\pi}^\pi \int_{-\pi}^\pi f_2(\lambda_1, \lambda_2)g(\underline{\lambda})d\lambda_1 d\lambda_2
    %  \end{align*}
Our estimator of this polyspectral mean is
\[
\widehat{M_g(f_2)} = 
T^{-3} \sum_{ |\tilde{\lambda}_1| \leq \mu_1, | \tilde{\lambda}_2| \leq \mu_2  } d( \tilde{\lambda}_1 )
        d( \tilde{\lambda}_2 )d(-\tilde{\lambda}_1 - \tilde{\lambda}_2 ).
%      g(\tilde{\lambda}_1, \tilde{\lambda}_2 ).
        % &= (2 \pi)^{-2} T^{-1} \int_{-\pi}^\pi \int_{-\pi}^\pi \sum_{t_1, t_2, t_3} X_{t_1}X_{t_2}X_{t_3}e^{-\iota \lambda_1 t_1 - \iota \lambda_2 t_2 + \iota (\lambda_1 + \lambda_2)t_3} \\
        % &\quad \quad \quad \quad \quad \quad \quad \quad \quad \quad \quad  \mathbb{I}(\lambda_1 \in (-h_1,h_1), \lambda_2 \in (-h_2, h_2)) d\lambda_1 d\lambda_2 \\
 %       &= (2 \pi)^{-2} \int_{-\mu_1}^{\mu_1} \int_{-\mu_2}^{\mu_2} I_2(\lambda_1, \lambda_2) d\lambda_1 d\lambda_2
    \]
The above example can be extended to the bispectral content on any general region $A \times B$, which will provide us the bispectral content in that region. 
We can also consider non-rectangular regions; for example, one can use an annular region like $\ind\left(0.2 \leq \lambda_1^2 + \lambda_2^2 \leq 0.6 \right)$.
%, which would give us the bispectral content in the given annular region. 
Several variations of 
such polyspectral means can be found in the literature; see \cite{boashash1993application, saidi2015application, martis2013application}.

\end{example}

\begin{example}[Example 5:]
\label{exam:g4}

Consider the Bartlett-type function $g(\lambda_1, \lambda_2) = (\pi - \lvert\lambda_1\rvert)(\pi - \lvert \lambda_2 \rvert)$, for which the estimate of the polyspectral mean is
\[
\widehat{M_g(f)} = {(2 \pi)}^2
T^{-3}\sum_{t_1, t_2, t_3} X_{t_1} X_{t_2} X_{t_3} \sum_{\tilde{\lambda_1}, \tilde{\lambda_2} \neq 0}
(\pi - \lvert \tilde{\lambda_1} \rvert)e^{-\iota \tilde{\lambda_1}(t_3 - t_1)}(\pi - \lvert \tilde{\lambda_2 }\rvert)e^{-\iota \tilde{\lambda_2}(t_3 - t_2)}.
\]
%
  % \begin{align*}
  %       \widehat{M_g(f)} &= (2 \pi)^2 T^{-2} \sum_{\lambda_1, \lambda_2 \neq 0} T^{-1} d(\lambda_1)d(\lambda_2)d(-\lambda_1 - \lambda_2)g(\lambda_1, \lambda_2) \\
  %       &\sim T^{-1} \int_{-\pi}^{\pi} \int_{-\pi}^{\pi} \sum_{t_1, t_2, t_3} X_{t_1} X_{t_2} X_{t_3} e^{-\iota \lambda_1 t_1 - \iota \lambda_2 t_2 + \iota (\lambda_1 + \lambda_2)t_3} (\pi - \lvert \lambda_1 \rvert)(\pi - \lvert \lambda_2 \rvert) \\
  %       &= T^{-1}\sum_{t_1, t_2, t_3} X_{t_1} X_{t_2} X_{t_3} \int_{-\pi}^\pi \int_{-\pi}^\pi (\pi - \lvert \lambda_1 \rvert)e^{-\iota \lambda_1(t_3 - t_1)}(\pi - \lvert \lambda_2 \rvert)e^{-\iota \lambda_2(t_3 - t_2)} 
  %   \end{align*}
Letting $c_h = (2 - 2 \cos (\pi h))/h^2$ for
$h \neq 0$ and $c_0 = \pi^2$, 
%Using trivial calculations, 
it can be shown that
$\widehat{M_g(f)}$ is asymptotically
equivalent to
$ 
 T^{-1}\sum_{t_1, t_2, t_3} X_{t_1} X_{t_2} X_{t_3} c_{t_3 - t_1} c_{t_3 - t_2}.
$ 
This bispectral mean can be interpreted as the spectral content in a pyramidal region around the center and tapered away from $(0,0)$.
% \[
% \pi^2 T^{-1}\sum_{t_1, t_2, t_3} X_{t_1} X_{t_2} X_{t_3} \frac{\sin\left(\frac{\pi (t_3 - t_1)}{2} \right)}{\frac{\pi (t_3 - t_1)}{2}} \frac{\sin\left(\frac{\pi (t_3 - t_2)}{2} \right)}{\frac{\pi (t_3 - t_2)}{2}}.
% \]
    % Now,
    % \begin{align*}
    %     & \quad \quad  \int_{-\pi}^\pi (\pi - \lvert \lambda \rvert) e^{-\iota \lambda (t_3 - t_1)}d \lambda \\
    %     &= \pi^2 \int_{-1}^1 (1 - \lvert z \rvert)e^{-\iota \pi z(t_3 - t_1)}dz \\
    %     &= \pi^2 \int_{-1}^1 \Pi(z) * \Pi(z) e^{-\iota \pi z(t_3 - t_1)}dz \\
    %     &= \pi^2 \frac{sin\left(\frac{\pi (t_3 - t_1)}{2} \right)}{\frac{\pi (t_3 - t_1)}{2}} 
    % \end{align*}
    % where $\Pi(t) = \mathbb{I}(\lvert t \rvert \leq 0.5)$, and $*$ denote the convolution function. Hence, the final expression would be:
    % \begin{align*}
    %     \widehat{M_g(f)} &= \pi^2 T^{-1}\sum_{t_1, t_2, t_3} X_{t_1} X_{t_2} X_{t_3} \frac{sin\left(\frac{\pi (t_3 - t_1)}{2} \right)}{\frac{\pi (t_3 - t_1)}{2}} \frac{sin\left(\frac{\pi (t_3 - t_2)}{2} \right)}{\frac{\pi (t_3 - t_2)}{2}} 
    % \end{align*}
%    The terms are non-zero only when $t_3 = t_1 + (2N-1)$ and $t_3 = t_2 + 2N-1$ for some   integer $N$, if we consider the $t_i$'s to be ordered.
\end{example}

\section{Asymptotic Results for Polyspectral Means }
\label{sec:results}
\setcounter{equation}{0}
% This section contains our main asymptotic results.
\subsection{Asymptotic mean and variance}
%The spectral mean is a classical parameter of interest in the time series literature. For example, in \citet{brillinger2012asymptotic} it was shown that the sample spectral mean converges to a Gaussian distribution with
%center given by the corresponding spectral mean, and  asymptotic variance involving
%the trispectrum.  
%A particular case is the sample autocovariance function (Example 1
% \ref{exam:g1}
%), which  is asymptotically normal with expectation given by
%the corresponding autocovariance, and
% variance involving a trispectral mean. 
 In this section, we derive asympttic expressions for the mean and the variance of the 
  the $k^{th}$ polyspectral mean estimate proposed in Section \ref{sec:expl}.
 % is asymptotically normal  by demonstrating  that the statistic's cumulants of order higher than $2$ tend to zero asymptotically. We  also compute the asymptotic variance term, which  involves polyspectra of order $2k + 1$. 
%We   already assume that $g$ is a symmetric weight function, and that the polyspectra of order up to $2k+1$ are continuous in $\lambda$. Also, let us assume that $\{ X_t \}$ is generated from a process with finite $(2k+1)$th order moment. 
We will need the following assumption on $k^{th}$ order autocumulants, as stated in \citet{brillinger1965introduction}.
\begin{equation}
\label{assumptionA}
\mbox{
\textbf{Assumption A[k]:}}\quad  
    \sum_{v_1, \ldots, v_{k-1}=-\infty}^\infty \lvert v_j \, \gamma (v_1, \ldots, v_{k-1})\rvert < \infty, \quad \textrm{ for } j = 1, \ldots, k-1.
\end{equation}
Assumption A[k] is directly related to the smoothness of the $k^{th}$ order polyspectra.

\begin{proposition}
\label{prop1}
Suppose that for some $k \geq 1$,
 $\{ X_t \}$ is a ${(k+1)}^{th}$  order stationary
time series 
%with finite ${(2k+1)}^{th}$  order moment,
and suppose that  Assumption A[r] holds for all $2\leq r\leq k+1$. 
Let $g$ be a  weight function 
satisfying the symmetry condition
(\ref{eq:sym-g}) that has finite absolute integral 
over the k-torus.
%$(-\pi,\pi)$. 
Then the polyspectral mean estimator $\widehat{M_g(f_k)}$
computed from a sample of length $T$ satisfies
(as $T \rightarrow \infty$)
\begin{align}
\label{eq:mean-convergence}
    E\widehat{M_g(f_k)} &= M_g\left(f_k\right) + o(1).
\end{align}
The $o(1)$-remainder term is $O(T^{-1})$, 
provided
\begin{align}
\label{eq:g-cond}
   \int |g - g_T| =  O\left(T^{-1}\right),
\end{align}
where $g_T (\omega)= \sum_{\underline{\tilde{\lambda}}} 
g(\underline{\tilde{\lambda}}) \ind\big(\omega\in \underline{\tilde{\lambda}} + (0, 2\pi T^{-1}]^k
\big)$,~
$\omega\in (-\pi, \pi]^k$ is the discretized version of $g$.

\end{proposition}

Condition \eqref{eq:g-cond} is a mild 
smoothness condition on $g$, allowing
Lipschitz functions (e.g.,  Examples 1 and 2)
as well as many  discontinuous functions
(e.g., Example 4). 
Proposition \ref{prop1} shows that the estimator's expectation tends to the  polyspectral mean,  
%at a rate of $T^{-1}$,
%The proof is quite simple, and is given in Section \ref{proof:prop1}. 
and hence it is 
%This result shows that the estimate defined for polyspectral mean is 
asymptotically unbiased.
%, and hence a valid estimate for polyspectral mean. 

Next, we seek to compute the
variance (or second cumulant) of the polyspectral mean. 
Following \cite{brillinger2001time}, a complex normal 
random variable $Z$ has 
%a complex normal distribution with
variance $\Sigma$ if $\Cum (Z,\Bar{Z}) = \Sigma$.
Since we are ultimately interested in demonstrating
 asymptotic normality of the spectral means, we shall first investigate 
 the variance,  $\Cum (Z,\Bar{Z}) $, 
 with $Z= \widehat{M_g(f)}$. 
%
%Since the polyspectral mean  can take a complex value, we need to take cumulant for complex-valued processes here. We can consider the cumulant in different directions, viz. $\Cum(X,X), \Cum(X, \Bar{X})$ or $\Cum(\Bar{X}, \Bar{X})$. All the cases will have similar proofs, and we will only look at the circular cumulant \citet{comon2010handbook}  $\Cum(X, \Bar{X})$ in this paper.
%
To that end, we first introduce some
convenient abbreviations.  As before, let
$\underline{\lambda} = \left(\lambda_1, \ldots, \lambda_k  \right)'$ and 
$[\underline{\lambda}]=\sum_{j=1}^k \lambda_j$,
but we can concatenate the two, 
% $-[\underline{\lambda}]$ to the end, 
which will be 
 denoted by $\{ \underline{\lambda} \}$,
 i.e., $\{ \underline{\lambda} \}'
 = \left( \underline{\lambda}, -[\underline{\lambda}]  \right)$.  Next, the $k$-vector of DFTs
 corresponding to $( d(\lambda_1), \ldots, d(\lambda_k))$ is abbreviated by $d(\underline{\lambda})$; similarly,
 $d \{ \underline{\lambda} \}$ is the $k+1$-vector
  given by $d(\underline{\lambda})$ with
  $d (- [ \underline{\lambda}])$ appended to the
  end.  Finally, let $\prod \underline{v}$ denote
  the product of all the components of any
  vector $\underline{v}$.  Then 
from (\ref{eq:mean-convergence}) it follows
that
\begin{align}
\label{cumulant1}
    \Cum \left(\widehat{M_g(f)}, \overline{\widehat{M_g(f)}} \right) = E\left( \widehat{M_g(f)}\overline{\widehat{M_g(f)}} \right) - M_g\left( f \right) \overline{M_g\left( f \right)} 
    + o(1).
    %+ O \left( T^{-1} \right).
\end{align}
The main challenge here is to derive 
an asymptotic expression
for the cross-product 
moment term in \eqref{cumulant1}, which requires some nontrivial combinatorial arguments. To highlight the key points we outline the main steps below,
 providing  details 
of the arguments  in Section \ref{proof:prop2}. 

Using our product notation, the cross-product 
moment
$E\left(\widehat{M_g(f)}\overline{\widehat{M_g(f)}} \right) $ can be written 
%\begin{align*}
% &E\left(\widehat{M_g(f)}
% \overline{\widehat{M_g(f)}} \right) \\
%    &= E\Bigg(T^{-2k}T^{-2}
%    \sum_{\tilde{\underline{\lambda}}, \tilde{\underline{\omega}}} \left( 2 \pi \right)^{-2k} d\left( \tilde{\lambda_1} \right)\cdots d\left(\tilde{\lambda_k} \right)d\left(- [\tilde{\underline{\lambda}}]\right) d\left(-\tilde{\omega_1}\right)\cdots d\left(-\tilde{\omega_k}\right)d\left( [\tilde{\underline{\omega}}]\right) \\
 %   & \quad \quad \quad \quad   \quad \quad \quad  g\left( \tilde{\underline{\lambda}} \right)
%    \overline{g(\tilde{\underline{\omega}})} \Phi( \tilde{\underline{\lambda}})
%    \Phi(\tilde{\underline{\omega}})\Bigg) \\
 \[
 T^{-2k-2}\left( 2 \pi \right)^{2k} \sum_{\tilde{\underline{\lambda}}, \tilde{\underline{\omega}}} g ( \tilde{\underline{\lambda}} ) \overline{g(\tilde{\underline{\omega}})}
    E\left( \prod d \{ \tilde{\underline{\lambda}} \}
     \prod d \{ \tilde{\underline{\omega}} \}  \right)
     \Phi(\tilde{\underline{\lambda}})
     \Phi(\tilde{\underline{\omega}}).
\]
%\end{align*}
Considering only the part inside the expectation in the above equation, one can show that for any $\la_1, \ldots, \la_k$ and $\omega_1,\ldots, \omega_k$, 
\begin{equation}
\label{eq:second}
    E\left( d\left(\tilde{\lambda_1}\right)\ldots d\left(\tilde{\lambda_k}\right)d\left(- [\tilde{\underline{\lambda}}]\right) d\left(-\tilde{\omega_1}\right)\ldots d\left(-\tilde{\omega_k}\right)d\left( [\tilde{\underline{\omega}}]\right) \right) 
    = \sum_{\sigma \in I_{2k+2}} \prod_{b \in \sigma} C_b,
\end{equation}
where $\sigma$ is a partition of $\{1, \ldots, 2k+2\}$
and $I_{2k+2}$ is the set of all partitions of $\{1, 2, \ldots, 2k+2\}$.  Each $C_b$ is a cumulant of DFTs whose
 frequencies $\lambda_j$ have subscript $j$ belonging to the set $b$    (which  is an element of $\sigma$).
From Lemma 1 of \citet{brillinger1967asymptotic}, under Assumption A[2k+1] (or (\ref{assumptionA})), each such term $C_b$ is non-zero if and only if the sum of the frequencies within the set is $0$. The presence of $\Phi(\lambda)$ and   $\Phi(\omega)$   ensures that the frequencies in the sum are  concentrated on a subset of
the $k$-dimensional torus with the property
that   $\sum_{j=1}^l \lambda_j \not \equiv 0 \textrm{ (mod }2\pi)$  and  $\sum_{j=1}^l \omega_j \not\equiv 0 \textrm{ (mod }2\pi)$ for all $1\leq l \leq k$. 
In other words, the Fourier frequencies do not lie in any sub-manifold, i.e., no subset of $\underline{\lambda}$ or $\underline{\omega}$ has their sum equal to $0$. Hence, the terms will be non-zero in two possible ways depending on partitions of the set
 $\{\{\underline{\lambda}\}, \{-\underline{\omega}\}\}$:
\begin{itemize}

\item The partition consists of the sets 
   $\{\underline{\lambda}\}$ and $\{-\underline{\omega}\}$, 
    \item The partition consists of sets, each
    of which has at least one element of
    both  $\{\underline{\lambda}\}$ and $\{-\underline{\omega}\}$. 
\end{itemize}
The first case yields an expression that
equals the second term of (\ref{cumulant1}), 
resulting in its cancellation out of the 
cumulant; so the only remaining terms will be those arising from the second case above. Now suppose we have $m$ such mixture partitions of $\{\underline{\lambda}\} = (\lambda_1, \ldots, \lambda_k, - [\underline{\lambda}])$ and
$\{ -\underline{\omega}\} = (-\omega_1, \ldots, -\omega_k,  [\underline{\omega}])$.
 We shall designate the sets of a partition
 through binary matrices, where each
 row corresponds to a set of the partition,
 and each column has a one or zero depending
 on whether the corresponding entry of
 $\{ \underline{\lambda} \}$ belongs
 to that set of the partition.  More 
 specifically, let  $\mathcal{B}_{m,k+1}$
equal  the set of $m \times (k+1)$-dimensional
matrices with binary entries; then for
$A, B \in \mathcal{B}_{m,k+1}$,
$A_{ij}$ (for $1 \leq i \leq m$ and 
$1 \leq j \leq k+1$) equals $1$ if
the $j^{th}$ component of  $\{ \underline{\lambda} \}$
is in set $i$ of the partition
(and is zero otherwise).  Likewise
$B_{ij}$ encodes whether the $j^{th} $ component of  $\{ \underline{\omega} \}$
is in the same set $i$ of the partition.
Because all the sets of a partition must
be disjoint, each column of such a matrix
$A$ (or $B$) has exactly one $1$; likewise,
 because we are considering partitions
 that have at least one element of both
  $\{\underline{\lambda}\}$ and $\{-\underline{\omega}\}$, it follows that
  each row of $A$ (and of $B$) must have
  at least one $1$.  Let us denote
 this subset of $\mathcal{B}_{m,k+1}$
 as $\widetilde{\mathcal{B}}_{m,k+1}$:
\[
\widetilde{\mathcal{B}}_{m,k+1} = 
\left\{ A \in \mathcal{B}_{m,k+1} \vert
 \textrm{ every row is non-zero and  every column has exactly one 1} \right\}.
\]
 Now the frequencies of a set in a partition
have the property that they sum to zero, and
this is the defining property of such
partitions; this essentially refines our
class of potential matrices 
$A, B \in \widetilde{\mathcal{B}}_{m,k+1}$.
 This constraint is simply expressed
 by $A \{\underline{\lambda}\} = B \{\underline{\omega}\} $.
 So it follows that the collection of
 partitions corresponding to
 given $\{ \underline{\lambda} \}$ and
 $\{ -\underline{\omega} \}$ is encoded by
\begin{equation}
\label{matCond}
\left\{  A, B \in \widetilde{\mathcal{B}}_{m,k+1}
\vert   A \{\underline{\lambda}\} = B \{\underline{\omega}\} \right\}.
\end{equation}
%
% \begin{align}
% \label{matCond}
% \nonumber \{ A \{\underline{\lambda}\} - B \{\underline{\omega}\} &= 0_{m \times 1} \bigg\vert 
%  A, B \in \mathcal{B}_{m,k+1}, 
% \textrm{  such that} \\
% & \textrm{ every row is non-zero and} \\
% \nonumber & \textrm{ every column has exactly one 1} \}.  
% \end{align}
% where $\{\underline{\lambda}\} = \{\lambda_1, \ldots, \lambda_k,  - [\underline{\lambda}]\}$, and $\{\underline{\omega}\} = \{\omega_1, \ldots, \omega_k,  - [\underline{\omega}]\}$.
As an illustration, for bispectra $k=2$ and $m=2$  the  set of possible choices of matrices $A$ and $B$ in
 (\ref{matCond}) is
\[\begin{pmatrix}
1 & 1 & 0\\
0 & 0 & 1
\end{pmatrix}, \begin{pmatrix}
1 & 0 & 1\\
0 & 1 & 0
\end{pmatrix}, \begin{pmatrix}
0 & 1 & 1\\
1 & 0 & 0
\end{pmatrix}, \begin{pmatrix}
1 & 0 & 0\\
0 & 1 & 1
\end{pmatrix}, \begin{pmatrix}
0 & 1 & 0\\
1 & 0 & 1
\end{pmatrix}, \begin{pmatrix}
0 & 0 & 1\\
1 & 1 & 0
\end{pmatrix}.
\]
For example, with the choices
\[
   A = \begin{pmatrix}
1 & 1 & 0\\
0 & 0 & 1
\end{pmatrix} \quad \mbox{and} \quad 
B = \begin{pmatrix}
1 & 1 & 0\\
0 & 0 & 1
\end{pmatrix},
\]
 we see that (\ref{matCond}) indicates two constraints, the first of which (corresponding
 to the upper row of the matrices) is
 $\lambda_1 + \lambda_2 = \omega_1 + \omega_2$;
 the second constraint from the lower row is
 identical to  the first due to the fact that the vectors $\{\underline{\lambda}\}$ and $\{\underline{\omega}\}$ add up to $0$.
 As a result, the corresponding partition
 consists of the sets  $(\lambda_1, \lambda_2, -\omega_1, -\omega_2)$ and $(-\lambda_1 - \lambda_2, \omega_1 +\omega_2)$, with a nonzero
 value of $\prod_{b\in \si} C_b$ provided  that
 $\lambda_1 + \lambda_2 = \omega_1 + \omega_2$. 

 In general we always have constraints given
 by the fact  that the sum of all the elements of $\{\underline{\lambda}\}$ and $\{\underline{\omega}\}$ are zero, since we have 
 $-[{\underline\lambda}]$
and $-[{\underline{\omega}}]$ 
 %%%%
 % $-\ceil{\underline\lambda}$ and $-\ceil{\underline{\omega}}$ 
 as the final
 components of $\{\underline{\lambda}\}$ and
 $\{\underline{\omega}\}$. 
 The rows of the binary matrices in 
 (\ref{matCond}) can be characterized 
 as follows:  let $p$ denote the   number of   possible ways to partition $k+1$ into positive integers, and define  $L_m$ to be the collection of all possible sets of positive integers that sums up to $k+1$.  That is,  
  $L_m = \{\underline{l}_1, \ldots, \underline{l}_p \}$,  where each $\underline{l}_j$ (for
  $1 \leq j \leq p$) is a vector of
 $m$ positive  integers whose components 
 sum to $k+1$.  (For instance, with
     $m=2$ and $k=2$, $L_m = \{(1,2), (2,1)\}$.)
 Then for every $\underline{l} = \{l (1), \ldots, l (m) \} \in L_m$, let
\[
    \zeta_{\underline{l}} = \left\{  A \in
\widetilde{\mathcal{B}}_{m,k+1} \vert 
    \textrm{ the $i^{th}$  row of $A$   has $l (i)$ ones }
 \right\},
\]
% \begin{equation*}
%     \Xi_{L_m, L_n} = \left\{ \underline{\lambda}, \underline{\omega} \Bigg\vert A \underline{\lambda} + B \underline{\omega} = 0 \textrm{ for some } A \in \zeta_{L_m}, B \in \zeta_{L_m'} \right\}  
% \end{equation*}
%
 so that $\cup_{\underline{l} \in L_m} \zeta_{\underline{l}} = \widetilde{\mathcal{B}}_{m,k+1} $.
 Let us define $f_\alpha(\underline{\lambda}_\alpha)$ to be the polyspectra of order $\alpha$ (where $\underline{\lambda}_\alpha$ is a vector of length $\alpha$), and let $r_{Aj}$ and $r_{Bj}$ denote the sum of the $j^{th}$ row of $A$ and $B$ respectively.  Further, set $r_j = r_{Aj} + r_{Bj} -1$, 
which denotes the sum of the $j^{th}$ row of $A$ and $B$  subtracted by 1.  Also let $\lambda_{r_{A_j}}$ denote the subset of $\{\underline{\lambda}\}$ that contains the elements corresponding to the non-zero positions of the $j^{th}$ row of $A$.
 Similarly, define $\omega_{r_{B_j}}$.
 Using the matrices from (\ref{matCond}), and some combinatorial arguments, it can be shown that the asymptotic variance $V$ of the polyspectral mean estimate is given by
 \begin{equation}
 \label{eq:asymp-var}
 V =  \sum_{m=1}^{k+1}\sum_{ \underline{l}, \underline{h} \in L_m} \sum_{A \in \zeta_{\underline{l}}, B\in \zeta_{\underline{h} } }
     \underbrace{\int_{-\pi}^\pi \ldots \int_{-\pi}^\pi}_{A\{\underline{\lambda}\} = B \{\underline{\omega}\} } g(\underline{\lambda})\overline{g(\underline{\omega})} \prod_{j=1}^m  \tilde{f}_{r_j } \left(\lambda_{r_{A_j}} , \omega_{r_{B_j}} \right)d\underline{\lambda}d\underline{\omega},
\end{equation}    
where  $\tilde{f}_k(\{\underline{\lambda}\}) = f_k(\underline{\lambda})$ is the $k^{th}$ order polyspectra. It is to be noted that the case when all the $\{\underline{\lambda}\}$ and $\{\underline{\omega}\}$ lie in the same partition is covered in the case $m=1$, and is simply the polyspectral mean of order $2k+1$ with weight function $g(\underline{\lambda})g(\underline{\omega})$, as mentioned in the introduction.  These derivations
    are summarized in the following proposition,
    whose formal proof is provided   in Section \ref{proof:prop2}.

\begin{proposition}
Suppose that, for some $k \geq 1$,
 $\{ X_t \}$ is a ${(2k+2)}^{th}$ order stationary
time series with finite ${(2k+2)}^{th}$ order moment,
and suppose that  Assumption A[r] holds for all $r=2,\ldots, 2k+2$. 
Let  $g$ be a  weight function 
satisfying conditions
(\ref{eq:sym-g}) 
and \eqref{eq:g-cond}, 
and let $g$  have 
a finite absolute integral 
over the k-torus.  Then the variance of the 
polyspectral mean estimator $\widehat{M_g(f_k)}$
based on a sample of length $T$ 
% Suppose $X_1, \ldots, X_T$ satisfies the zero-mean stationarity conditions, and the weight function g satisfies the symmetry condition. Furthermore, suppose Assumption 1 is satisfied. Then, the second cumulant 
is of the form $T^{-1}V + O(T^{-2})$, where
$V$ is given by (\ref{eq:asymp-var}).
% \begin{align*}
%    V =  \sum_{m=1}^{k+1}\sum_{ \underline{l}_m , \underline{l}_m' \in L_m} \sum_{A \in \zeta_{l_m}, B\in \zeta_{l_m'}} 
%     % \tau_{\underline{l}_m}\tau_{\underline{l}_m'}
%     \underbrace{\int_{-\pi}^\pi \ldots \int_{-\pi}^\pi}_{A\{\underline{\lambda}\} - B \{\underline{\omega}\} = 0} g(\underline{\lambda})\overline{g(\underline{\omega})} \prod_{j=1}^m  \tilde{f}_{r_j } \left(\lambda_{r_{A_j}} , \omega_{r_{B_j}} \right)d\underline{\lambda}d\underline{\omega}
% \end{align*}
\end{proposition}

%\begin{remark} 
%   The case $m=1$ in the variance formula
%   (\ref{eq:asymp-var})   is given by
%\[
%\int_{-\pi}^\pi \ldots \int_{-\pi}^\pi g(\underline{\lambda})\overline{g(\underline{\omega})}f_{2k+1}
%(\lambda_1, \ldots, \lambda_k, - [\underline{\lambda}], \omega_1, \ldots, \omega_k)
%d\underline{\lambda}d\underline{\omega}.
%\]
%\end{remark}

%Therefore, we have obtained the asymptotic distribution of the estimated polyspectral mean for the general weight function.
A simple extension of the proof of Proposition 2
shows that the covariance between polyspectral means of different orders goes to $0$ at rate  $T^{-1}$. The proof is given in Section \ref{proof:corr1}.

\begin{corro}
\label{cor:theo1} 
Suppose that, for some $k \geq 1$,
 $\{ X_t \}$ is a ${(2k+3)}^{th}$ order stationary
time series with finite ${(2k+3)}^{th}$ order moment,
and suppose that  A[r] holds for all $r=2,\ldots, 2k+3$.
Let $g_1$ and $g_2$ be weight functions that each
satisfy the symmetry condition
(\ref{eq:sym-g}), and have
finite absolute integral 
over the k- and (k+1)-toruses respectively.  Then
$T \mbox{Cov} (\widehat{M_{g_1}(f_k)}, \overline{\widehat{M_{g_2}(f_{k+1})}})$ 
tends to 
\[
 \sum_{m=1}^{k+1}
 \sum_{ \underline{l} \in L_m^{(k)}, \underline{h} \in L_m^{(k+1)}}
 \sum_{A \in \zeta_{\underline{l}}, B\in \zeta_{\underline{h}}} 
    \underbrace{\int_{-\pi}^\pi \ldots \int_{-\pi}^\pi}_{A\{\underline{\lambda}\} =
    B \{\underline{\omega}\} } 
g_{1}\left(\underline{\lambda}_k\right)\overline{g_2(\underline{\omega}_{k+1})}\prod_{j=1}^m f_{r_j }\left(\lambda_{r_{A_j}} , \omega_{r_{B_j}} \right)  d\underline{\lambda}d\underline{\omega},
\]
% {\bf This does not look right. $\lambda$ and $\omega$ have different
% orders - viz. $k+1$ and $k+2$. I think in the second sum, $l\in L_m^{(k)}$ and $h\in L_m^{(k+1)}$, with obvious interpretation of the superscripts. Need to define / make it precise!}
% \begin{align}
%     \nonumber Cov(\widehat{M_{g_1}(f_k)}, \overline{\widehat{M_{g_2}(f_{k+1})}}) =& T^{-1} \Bigg\{\sum_{m=1}^{k+1}\sum_{ \underline{l}_m , \underline{l}_m' \in L_m} \sum_{A \in \zeta_{l_m}, B\in \zeta_{l_m'}} 
%     % \tau_{\underline{l}_m}\tau_{\underline{l}_m'}
%     \underbrace{\int_{-\pi}^\pi \ldots \int_{-\pi}^\pi}_{A\{\underline{\lambda}\} - B \{\underline{\omega}\} = 0} \prod_{j=1}^m \\
%     & g_{1}\left(\underline{\lambda}_k\right)\overline{g_2(\underline{\omega}_{k+1})}\prod_{j=1}^m f_{r_j }\left(\lambda_{r_{A_j}} , \omega_{r_{B_j}} \right) \Bigg\}
% \end{align}
where $L_m^{(k)} = \{\underline{l}_1, \ldots, \underline{l}_p \}$,  where each $\underline{l}_j$ (for
  $1 \leq j \leq p^{(k)}$) is a vector of
 $m$ positive  integers with $\ell_1$ norm=$k+1$, and $p^{(k)}$  is 
 the   number of   possible ways to partition $k+1$ into positive integers.
\end{corro}

\subsection{Asymptotic normality}
We shall give two results on asymptotic normality of the polyspectral
mean estimators under two different sets of conditions. The first approach 
is based on the convergence of the sequence of cumulants 
(which is also known as the method of moments).  Now that we have derived the asymptotic  mean and variance of the polyspectral mean estimator,  we may  prove its asymptotic normality by showing that the higher order autocumulants (of order 3 and more) go to $0$ at a rate faster than $T^{-1}$. 
To do that we  make use of  Theorem 2.3.3 of \citet{brillinger2001time}, which states that
\beq
   \Cum\Big(\prod_{j=1}^{J_1} X_{1j}, \ldots, \prod_{j=1}^{J_I} 
        X_{Ij}\Big) = \sum_{\nu} \Cum(X_{ij}; ij \in \nu_1) \cdots 
   \Cum(X_{ij} ; ij \in \mu_p), 
\label{eq:Cum}
\eeq
where the summation is over all indecomposible partition $\nu = \nu_1 \cup \ldots \cup \nu_p$ of a 2-D array with $I$ rows, with the $i$th row 
having $J_i$ elements.  
%$$
%\centering
%\begin{pmatrix}
%(1,1) & \cdots & (1,J_1) \\
%\cdot & & \cdot \\
%\cdot & & \cdot \\
%\cdot & & \cdot \\
%(I,1) & \cdots & (I,J_I) 
%\end{pmatrix}
%$$ 
Using this result, it can be shown that 
cumulants of order $r\geq 3$ tend to $0$ at  rate $T^{-r+1}$ (see Section \ref{proof:theorem1} for a detailed proof), and
we have the following result:

\begin{theorem}
\label{theo1}
Suppose that
 $\{ X_t \}$ is a ${(k+1)}^{th}$ order stationary
time series 
%with finite 
%${(2k+1)}^{th}$ order 
%{moments of all orders}
and suppose that  Assumption A[$k$] holds,
for all $k\geq 1$. 
Let $g$ be a  weight function 
satisfying  conditions
(\ref{eq:sym-g}) and \eqref{eq:g-cond},
and let $g$  have a finite absolute integral over the k-torus. 
% Suppose $X_1, \ldots, X_T$ is a sample from a stationary time series, and $M_g(f)$ and $\widehat{M_g(f)}$ be as defined earlier. Let g be a Hermitian weight function such that $\int_{\underline{\lambda}} \lvert g\left( \underline{\lambda} \right)\rvert d\underline{\lambda} < \infty$. Let $A$ and $B$ are $m \times (k+1)$ binary matrices as defined in (\ref{matCond}).
% such that $A \{\underline{\lambda}\} - B \{\underline{\omega}\} = 0_{m \times 1}$, with the following constraints:
% \begin{itemize}
%     \item Rows are disjoint.
%     \item Every column must have exactly one 1.
% \end{itemize}
% Further let us define $f_\alpha(\underline{\lambda}_\alpha)$ to be the polyspectra of order $\alpha$ ($\underline{\lambda}_\alpha$ is a vector of length $\alpha$), $r_{Aj} $ and $r_{Bj} $ denote the sum of the $j^{th}$ row of $A$ and $B$ respectively, $r_j  = r_{Aj}  + r_{Bj} -1$, and $\lambda_{r_{A_j}}$ denotes the subset of $\{\underline{\lambda}\}$ which contains the elements corresponding to the non-zero positions of $r^{th}$ row of $A$. Similarly, define $\omega_{r_{B_j}}$. Let $L_m$ and $L_m'$ be the set of all possible rows of $A$ and  $B$.
    Then 
   \[
   \sqrt{T} \left(\widehat{M_g(f_k)} - M_g(f_k) \right)  \Rightarrow \mathcal{N}(0, V)
  \]
    as $T \rightarrow \infty$,     where 
$V$ is given by (\ref{eq:asymp-var}).
    %
    % \begin{align}
    % \label{varform}
    %  V =  \sum_{m=1}^{k+1}\sum_{ \underline{l}_m , \underline{l}_m' \in L_m} \sum_{A \in \zeta_{l_m}, B\in \zeta_{l_m'}} 
    % % \tau_{\underline{l}_m}\tau_{\underline{l}_m'}
    % \underbrace{\int_{-\pi}^\pi \ldots \int_{-\pi}^\pi}_{A\{\underline{\lambda}\} - B \{\underline{\omega}\} = 0} g(\underline{\lambda})\overline{g(\underline{\omega})} \prod_{j=1}^m  \tilde{f}_{r_j } \left(\lambda_{r_{A_j}} , \omega_{r_{B_j}} \right)d\underline{\lambda}d\underline{\omega}
    % \end{align} 
\end{theorem}

Thus, the polyspectral mean estimators are asymptotically normal with  asymptotic variance that involves the $(2k+1)$th order polyspectral
density. However, one potential drawback of this approach is that it requires all moments of $X_t$ to be finite, which is too strong for many applications. It is possible to relax the moment condition by
using a mixing condition and/or using the approach of 
\citet{wu2005pnas}
based on the physical dependence 
measure.  
For the second result on asymptotic normality,
here we follow neither of these approaches.
We develop a new approach 
based on a suitable structural condition on $\{X_t\}$ 
(similar to \citet{wu2005pnas}, but with a different
line of proof) 
%of proving asymptotic normality 
that is valid under a weaker moment condition. 
In fact, the moment condition is minimal, as it 
requires 
existence of $(2k+2)$th order moments only,
which  appears in the asymptotic variance of 
$\sqrt{T} \big(\widehat{M_g(f)} - M_g(f) \big)   $
(cf. Proposition 2)  and must necessarily 
be finite. The 
structural conditions on the sequence $\{X_t\}$ 
allow us to
approximate the polyspectral mean estimators by its analog based on a 
suitably constructed approximating process that is simpler, 
and that has all finite moments with 
quickly decaying  cumulants of all orders. To describe the 
set up, let 
$\{\ep_i\}_{i\in \bbz}$ be a collection of independent 
and identically distributed random variables, with $\bbz=\{0, \pm 1,\ldots \}$, 
and let 
\beq
X_t= a( \ldots, \ep_{t-1}, \ep_t)
\label{eq:nonlin}
\eeq
 for some Borel measurable function $a:\bbr^\infty \raw \bbr$. The class of 
 stochastic processes given by \eqref{eq:nonlin} includes common linear
 time series models such as ARMA ($p, q$), as well as many nonlinear 
 processes, such as the threshold AR process and the ARCH process,
 among others. See \citet{HsingWu2004} for a long list of 
 processes that can be represented as \eqref{eq:nonlin}. 

% For $m\geq 1$, set 
% $X_{t,m} = a(\ldots, 0,0,\ep_{t-m},\ldots, \ep_t)$
% and define 
% %\beq
%   $   a_r(m) = [E|X_1 - X_{1,m}|^{r}]^{1/r}
% $, $r\geq 1$.   
% %\label{eq:g-approx}
% %\eeq
% We shall suppose that $a_r(m)$ is nonincreasing in $m$.
% %\footnote{What is $c(m)$, do you mean $a_r(m)$?
% %Is it nonincreasing in $m$, rather than $n$?} 
%  The following result 
% proves the asymptotic normality of the polyspectral 
% mean estimators under a decay condition on 
% $a_r(m)$. Assumptions A[k]-s are not needed. 

% \begin{theorem}
% \label{theo1p}
% Suppose that $\{X_t\}$ admits the representation in 
% \eqref{eq:nonlin} and that $E|X_1|^{2k+2}<\infty$
% for some $k\geq 1$.  
% Let $g$ be a  weight function 
% satisfying  conditions
% (\ref{eq:sym-g}) and \eqref{eq:g-cond},
% \textcolor{red}{and let} $g$ have a finite absolute integral over the k-torus. If, in addition, with $r=2k+2$,
% \beq
% \sum_{m\geq 1} a_r(m) m^{r} <\infty,
% \label{eq: c-cond}
% \eeq
% then 
%  $ 
%    \sqrt{T} \left(\widehat{M_g(f)} - M_g(f) \right)  \Rightarrow \mathcal{N}(0, V)
%   $ 
%     as $T \rightarrow \infty$,     where 
% $V$ is given by (\ref{eq:asymp-var}).
% \end{theorem}

The next corollary gives the joint distribution of multiple polyspectral mean estimators of the same order but with different weight functions.
%It can be shown that the polyspectral mean estimates of different weight functions asymptotically follow a multivariate normal distribution with a non-diagonal covariance matrix, given by (\ref{covTerm}).
The proof is discussed in the first part of Section \ref{theorem2Proof}. Corollary \ref{covCorr} will be useful in deriving the time series linearity test  proposed  in the next section.

% The resulting form of the asymptotic covariance will be:

% \begin{equation*}
%     T^{-1} \left\{\sum_{m=2}^{k+1}\sum_{\underline{l_m} \in L_m} \tau_{\underline{l}_m}\tau_{\underline{l}_m'} T^{-2k + m -2}\sum_{\underline{\lambda},\underline{\omega} \in \Xi_{L_m, L_n}}g_{1}\left(\underline{\lambda}_k\right)g_{2}\left(\underline{\omega}_{k+1}\right)\prod_{j=1}^m f_{r_j }\left(\lambda_{r_{A_j}} , \omega_{r_{B_j}} \right) \right\}
% \end{equation*}

% The corollary can be stated as:
% \begin{corro}
% \begin{align}
%     \nonumber Cov(M_{g_1}(\hat{f}_k), \overline{\widehat{M_{g_1}(f_k)}) =& T^{-1} \Bigg\{\sum_{m=2}^{k+1}\sum_{\underline{l_m} \in L_m} \tau_{\underline{l}_m}\tau_{\underline{l}_m'} T^{-2k + m -2} \\
%     &\sum_{\underline{\lambda},\underline{\omega} \in \Xi_{L_m, L_n}}g_{1}\left(\underline{\lambda}_k\right)g_{2}\left(\underline{\omega}_{k+1}\right)\prod_{j=1}^m f_{r_j }\left(\lambda_{r_{A_j}} , \omega_{r_{B_j}} \right) \Bigg\}
% \end{align}
% \end{corro}

\vspace{-0.05 cm}
\begin{corro}
\label{covCorr}
Suppose the assumptions of either Theorem \ref{theo1}
or Theorem \ref{theo1p}  hold, and  consider the polyspectral means of order $k$ for different weight functions $g_1, \ldots, g_d$.  Then, the $d$-dimensional vector 
\[
\sqrt{T} \left( \widehat{M_{g_1}(f_k)} - M_{g_1}(f_k), \ldots, \widehat{M_{g_d}(f_k)} - M_{g_d}(f_k)
\right)
\]
  converges to a $d$-variate normal distribution with mean 0 and covariance matrix with entry $i,n$ given by  $V_{g_i, g_n}$
  where 
\begin{equation}
\label{covTerm}
V_{g_i,g_n} =  \sum_{m=1}^{k+1}\sum_{ \underline{l}, \underline{h} \in L_m} 
     \sum_{A \in \zeta_{\underline{l}}, 
     B\in \zeta_{\underline{h}}} 
     \underbrace{\int_{-\pi}^\pi \ldots \int_{-\pi}^\pi}_{A\{\underline{\lambda}\} = B \{\underline{\omega}\} } g_i(\underline{\lambda})
     \overline{g_n (\underline{\omega})}  \prod_{j=1}^m  \tilde{f}_{r_j } \left(\lambda_{r_{A_j}} , \omega_{r_{B_j}} \right)d\underline{\lambda}d\underline{\omega}.
\end{equation}
\end{corro}

\section{Testing of Linear Process Hypothesis Using Bispectrum 
\label{sec:application}}
\setcounter{equation}{0}
% \vspace{-0.5cm}
As discussed earlier, it is often of interest to determine whether a process is linear. \citet{ghosh2021} showed that it is possible to use quadratic prediction instead of the widely used linear predictor to get significant improvement in prediction error. In this section we will provide a novel test of linearity using the  bispectrum. 
Consider the null hypothesis
\begin{equation}
    \label{eq:null-linear}
    H_0 : \{X_t\} \, \mbox{is a linear process of the form} \, X_t = \psi(B)\epsilon_t,
\end{equation}
where $B$ is the backshift operator.
 Here, the moving average filter $\psi (B) =
 \sum_j \psi_j B^j$ is assumed to be known
 under $H_0$.  Thus, when the null hypothesis is
 true,  the bispectrum will be of the form $f_2(\lambda, \omega) = \mu_3 \psi(e^{-\iota \lambda})\psi(e^{-\iota \omega}) \psi(e^{\iota (\lambda + \omega)})$ for some $\mu_3$.  Setting
 $\Psi (\lambda, \omega) = 
 \psi(e^{-\iota \lambda})\psi(e^{-\iota \omega}) \psi(e^{\iota (\lambda + \omega)})$,
 % $\Psi (\lambda, \omega) = f_2 (\lambda, \omega)/\mu_3$,
 %we construct a statistic 
we study the quantity:
\begin{equation*}
    \mathcal{T}(\lambda, \omega) = \frac{f_2(\lambda, \omega)}{\psi(e^{-\iota \lambda})\psi(e^{-\iota \omega}) \psi(e^{\iota (\lambda + \omega)})} =\frac{f_2(\lambda, \omega)}{\Psi(\lambda, \omega)},
\end{equation*}
noting that $\mathcal{T}(\lambda, \omega)$ will be constant under  (\ref{eq:null-linear}). Setting $g_{j,k}(x_1, x_2) = 
\exp \{ \iota j x_1 + \iota k x_2 \} /\Psi(x_1, x_2)$, it follows that 
for any $j$ and $k$ the expression
\begin{equation}
\label{eq:linparam-def}
\int_{-\pi}^\pi \int_{-\pi}^\pi \mathcal{T}(\lambda, \omega) e^{\iota \lambda j}e^{\iota \omega k} \,  d\lambda \, 
d \omega
= M_{g_{j,k}} (f_2)
%\langle \langle \mathcal{T} \rangle_j \rangle_k
\end{equation}
equals $0$ whenever either $j$ or $k$ (or both) 
are non-zero.  We will examine the case where
(\ref{eq:linparam-def}) is zero for a
 range of $j,k$ between $0$ and $M$
 (an upper threshold chosen by the
 practitioner) such that not both are zero;
 this set is described via
$ 
\mathcal{v} = \{(0,1),(0,2),\ldots(0,M),(1,0),\ldots(1,M),\ldots(M,0),\ldots(M,M)\}
$. 
Substituting  \tr{$\hat{f}_2$}  for \tr{$f_2$} in
(\ref{eq:linparam-def}),  we can construct a   test statistic of the form
%\begin{align*}
  \tr{$  \sum_{(j,k)\in \mathcal{v}}  T \big\lvert 
%    \langle \langle \hat{\mathcal{T}} \rangle_j \rangle_{k} 
\widehat{M_{g_{j,k}} (f_2)}
\big\rvert ^2$}.
%\end{align*}
This test statistic, after scaling by $V_{g_{j,k}, g_{j,k}}$ 
% (defined in (\ref{covTerm})
, i.e., 
the asymptotic variance corresponding to the weight function $g_{j,k}$), becomes
\begin{equation}
\label{eq:blt-stat}
    \mathcal{T}_{BLT} = \sum_{(j,k) \in \mathcal{v}}  \frac{T\lvert \widehat{M_{g_{j,k}}(f_2)}\rvert^2}{ V_{g_{j,k}, g_{j,k}}}.
\end{equation}
Note that under $H_0$, $V_{g_{j,k}, g_{j,k}}$ is known
 and given by (\ref{eq:asymp-var}),
 and hence  $\mathcal{T}_{BLT}$ is computable.
 As discussed in Corollary \ref{covCorr}, the asymptotic distribution of a vector of polyspectral means with different choices of $g$ function is a multivariate normal with covariance term given by (\ref{covTerm}). %Hence, suppose we consider the functions  $g_{j,k}(x_1, x_2) = 
%\exp \{ \iota j x_1 + \iota k x_2 \} /\Psi(x_1, x_2)$, where $j$ and $k$ runs from $0$ to $M$, such that both are not zero. 
%
In fact, under the null hypothesis the vector of scaled polyspectral means \tr{$\big\{ 
 T^{1/2} \widehat{M_{g_{\mathbf{v}}}} \vert \mathbf{v} \in \mathcal{v} \big\}$} converges in distribution to a multivariate normal variable with known covariance matrix $\mathbb{CV}_{BLT}$ defined by
  (\ref{cvblt}) below. 
 
 There are many applications where researchers have used linear processes to model time series data, and our proposed test can be used to identify deviations from the linear model. For example in Section \ref{sec:data} we consider  Sunspot data to which \citet{abdel2018statistical} fitted an AR(1) model; our proposed test can be used to   check whether the fitted model is 
 appropriate for the dataset, or whether
 there are departures from the bispectrum
 implied by such a model.  
% We can test for linearity by a two step process:
% \begin{itemize}
%     \item First fit an ARMA model to the data and obtain an estimate of the underlying linear process, say $\tilde{\psi}(\cdot)$. 
%     \item Apply the proposed linearity testing to test $H_0: X_t = \tilde{\psi}(B)\epsilon_t$
% \end{itemize}
%
Our test considers $\psi(\cdot)$ to be known, and as such the limiting distribution under $H_0$ is completely known. Then we have the following:
% \vspace{-0.5cm}
%%%%%%%%%%%%%%%%%%%%%%%%%%%%%%%%%%%%%%%%%%%%%%%%%%%%%%%%%%%%%%%%%%%%
\begin{theorem}
\label{thm:lin-test}
Let $\{X_t\}$ be a 
%third-order stationary
time series 
satisfying
%Assumption [A3]. 
%Then, under 

the null
hypothesis $H_0: X_t = \psi(B)\epsilon_t$ for 
i.i.d. random variables $\{\ep_t\}$ 
 with $E\ep_1=0$ and  $E\epsilon_1^6 <\infty$
and for  some known $\psi$
satisfying
$\sum_j |j|^7 |\psi_j|<\infty$. Then, 
the test statistic
defined in (\ref{eq:blt-stat}) is
asymptotically distributed as
$\sum_{j=1}^{(M+1)^2-1} \nu_j \zeta_j$, where $\zeta_j$ are i.i.d. $\chi^2_1$ random variables, and $\nu_1, \ldots, \nu_{(M+1)^2-1}$ are the eigenvalues  of 
the covariance matrix $\mathbb{CV}_{BLT}$,
with $(a,b)$th  entry 
\begin{equation}
\label{cvblt}
  \mathbb{CV}_{BLT}(a,b) \equiv  \mathbb{CV}_{BLT}^{(M)}(a,b) =
    \begin{cases}
      1 & \text{ if } a=b\\
      \frac{ V_{g_{j_a,k_a}, g_{j_b,k_b}} }{\sqrt{V_{g_{j_a,k_a}, g_{j_a,k_a}}       
       V_{g_{j_b,k_b}, g_{j_b,k_b}} } } & \text{if } a \neq b.
    \end{cases}       
\end{equation}  
\end{theorem}
%%%%%%%%%%%%%%%%%%%%%%%%%%%%%%%%%%%%%%%%%%%%%%%%%%%%%%%%%%%%%%%%

\section{Simulation}
\label{sec:sim}

The simulation processes, weighting
functions for polyspectral means,
and simulation results are  
described below.

\subsection{AR(2) Process}

Our first simulation process is an AR(2) process $\{ X_t \}$ that is similar to Example 4.7 of \citet{shumway2005time}, and is  defined via 
$X_t = X_{t-1} - 0.9 X_{t-2} + \epsilon_t$, where $\epsilon_t \sim \mbox{Exp}(1)-1$.
%($r =1.1, n=100$).
In this case, $\phi(z) = 1 - z + 0.9 z^2$ and $\theta(z) =1$, where $\phi(B) X_t = \theta(B) \epsilon_t$.
%Also,
%\begin{align*}
%        \lvert \phi(e^{- 2\pi \iota \omega}) \rvert^2 = 2.81 - 3.8 cos(2\pi \omega) + 1.8 cos(4\pi \omega)
%\end{align*}
It follows that the spectral density is 
\begin{align*}
        f(\lambda) = \frac{1}{\phi(e^{-   \iota \lambda})\phi(e^{ \iota \lambda})} = \frac{1}{2.81 - 3.8 cos(\lambda) + 1.8 cos(2 \lambda)}.
\end{align*}
    Similarly, the bispectra and trispectra are  given by
    $ 
        f_2(\lambda, \omega) = {2}[{\phi(e^{-\iota \lambda})\phi(e^{-  \iota \omega})
        \phi(e^{ \iota (\lambda+\omega)})}]^{-1}$ and $
        f(\lambda_1,\lambda_2, \lambda_3) = {9}[{\phi(e^{-  \iota \lambda_1})
        \phi(e^{-  \iota \lambda_2})\phi(e^{-   \iota \lambda_3})
        \phi(e^{  \iota (\lambda_1 + \lambda_2 + \lambda_3)})}]^{-1}$,  respectively.
 %%%   
We also consider this process with
$\epsilon_t \sim \chi^2_4 - 4$, which changes
the scale of the polyspectra.

    %   The bispectra of the AR(2) process is depicted in Figure \ref{fig:bispecArma}.
    
    % \begin{figure}[h!]
    %     \centering
    %     \includegraphics[width = 0.45 \textwidth]{Figures/bispecARpersp.png}
    %     \includegraphics[width = 0.45\textwidth]{Figures/bispecARheatmap.png}
    %     \caption{Bispectra of the AR(2) Process}
    %     \label{fig:bispecArma}
    % \end{figure}
    
    % \item \textbf{Model 2 (AR(2)):} $X_t = X_{t-1} - 0.9 X_{t-2} + \epsilon_t$, where $\epsilon_t \sim \chi^2_4 - 4$ (r =1.1, n=100) (Similar to Example 4.7 of \citet{shumway2005time}).
    
    % In this case, $\phi(z) = 1 - z + 0.9 z^2$, $\theta(z) =1$, where $\phi(z) X_t = \theta(z) \epsilon_t$.
    
    % The spectra, bispectra and trispetra are as follows:
    
    % \begin{align*}
    %     f(\omega) &= \frac{8}{\phi(e^{- 2 \pi \iota \lambda})\phi(e^{2 \pi \iota \lambda})} \\
    %     f(\lambda, \omega) &= \frac{32}{\phi(e^{- 2 \pi \iota \lambda})\phi(e^{- 2 \pi \iota \omega})\phi(e^{2 \pi \iota (\lambda+\omega)})} \\
    %     f(\lambda_1,\lambda_2, \lambda_3) &= \frac{384}{\phi(e^{- 2 \pi \iota \lambda_1})\phi(e^{- 2 \pi \iota \lambda_2})\phi(e^{- 2 \pi \iota \lambda_3})\phi(e^{2 \pi \iota (\lambda_1 + \lambda_2 + \lambda_3)})}
    % \end{align*}

 % \vspace{-1.2 cm}

\subsection{ARMA(2,1) Process}
\label{sec:arma}
 \vspace{-0.4 cm}

    \begin{center}
   \begin{figure}[htbp]
        \centering
        \includegraphics[width = 0.45 \textwidth]{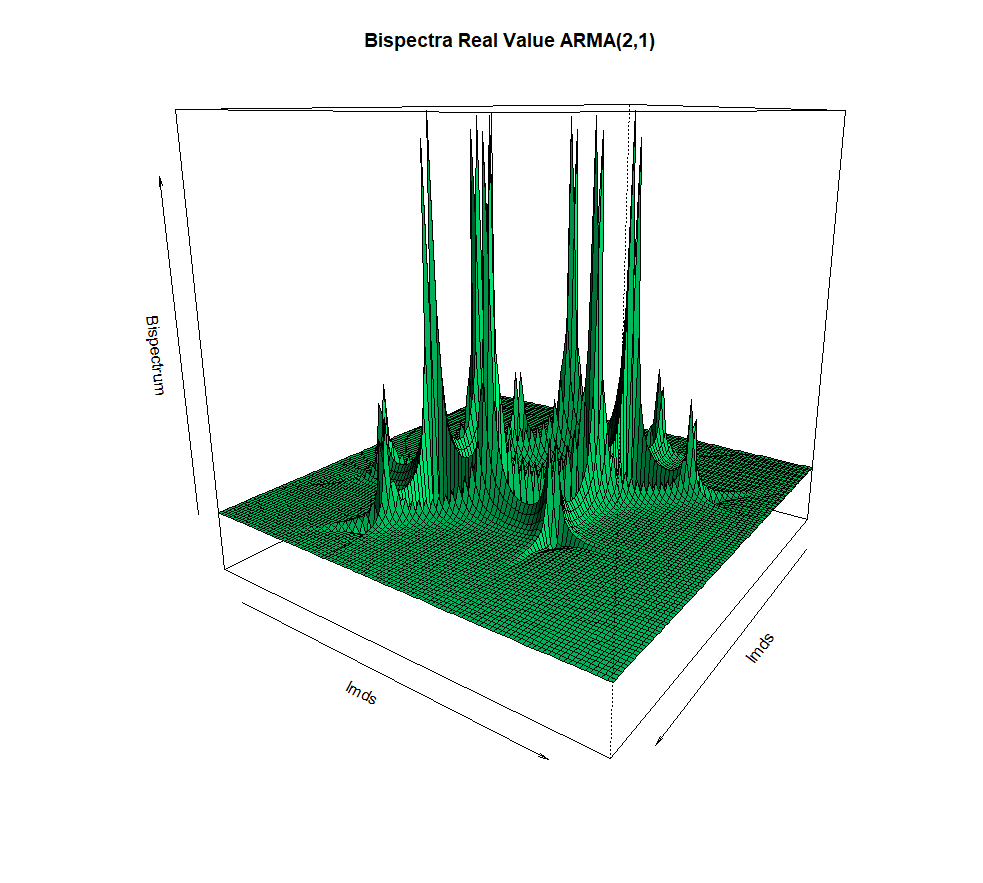}
        \includegraphics[width = 0.45\textwidth]{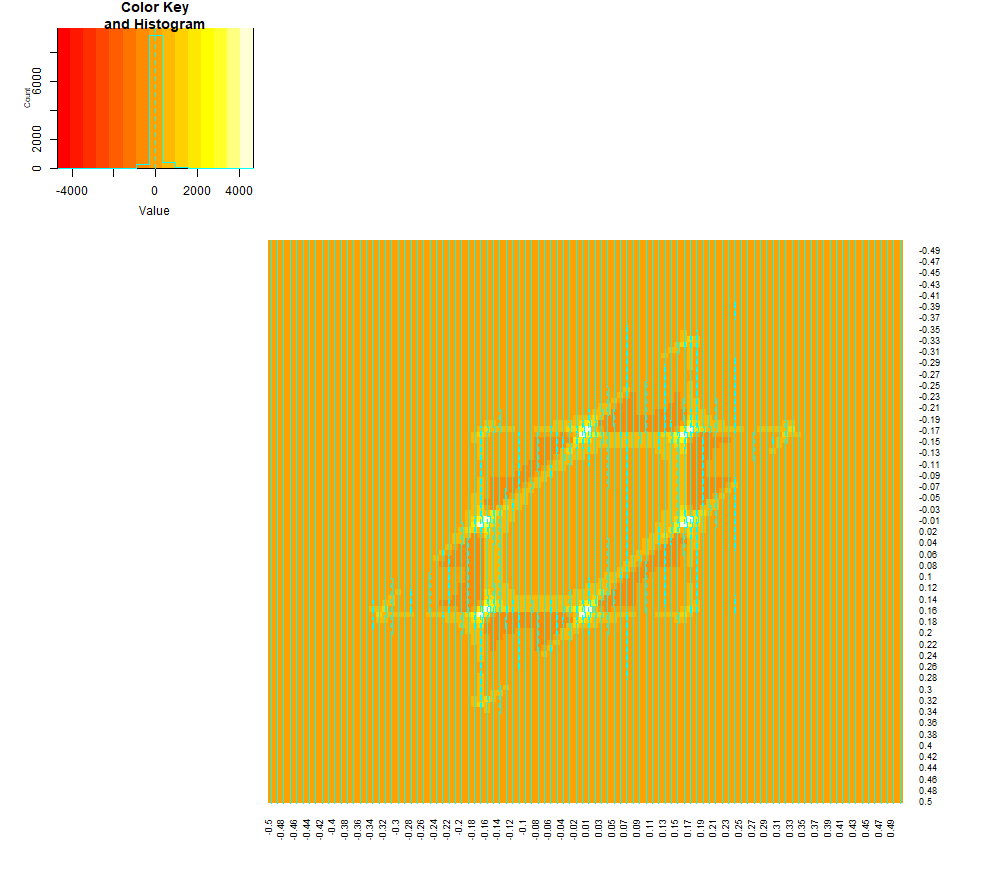}
        \caption{Bispectrum of ARMA(2,1) process with parameters given in Section \ref{sec:arma}, left panel exhibiting the 3d-plot and the right panel showing the heatmap.}
        \label{fig:bispecArma}
    \end{figure}     
    \end{center}
Next, we consider an ARMA(2,1) process $\{X_t\}$, defined via  $X_t -  X_{t-1} + 0.9 X_{t-2} = \epsilon_t + 0.8 \epsilon_{t-1}$, where $\epsilon_t \sim \mbox{Exp}(1)-1$, i.e. $\phi(z) = 1 - z + 0.9 z^2$ and $\theta(z) = 1 + 0.8 z$, where $\phi(B)X_t = \theta(B)\epsilon_t$. Hence, $\Psi(z) = \theta(z)/\phi(z) = (1 + 0.8z)/(1 - z + 0.9 z^2)$. The spectra, bispectra, and trispectra are:
    \begin{align*}
        &f(\lambda) = \lvert\Psi(e^{-   \iota \lambda})\rvert^2, \quad 
        f(\lambda, \omega) = 2\Psi(e^{-   \iota \lambda})\Psi(e^{-  \iota \omega})\Psi(e^{ \iota (\lambda+\omega)}), ~~\mbox{and} \\
        & f(\lambda_1,\lambda_2, \lambda_3) = 9\Psi(e^{-  \iota \lambda_1})\Psi(e^{-  \iota \lambda_2})\Psi(e^{-  \iota \lambda_3})\Psi(e^{  \iota (\lambda_1 + \lambda_2 + \lambda_3)}).
    \end{align*}
    The bispectra of the ARMA(2,1) process is depicted in Figure \ref{fig:bispecArma}.
    We also consider this process with
$\epsilon_t \sim \chi^2_4 - 4$, which changes
the scale of the polyspectra.

\subsection{Squared Hermite}

The final example is a squared Hermite process defined by $X_t = J_1 \, H_1 (Z_t) + J_2 \, H_2 (Z_t) = J_1 \, Z_t + J_2 \, Z_t^2 - J_2$, where $J_1 = 2$, $J_2 = 5$, and $\{ Z_t \}$ is a MA(1) process such that $Z_t = \epsilon_t + 0.4 \epsilon_{t-1}$ with $\epsilon_t \sim \mbox{Exp} (1) - 1$. The autocumulants are 
    {\allowdisplaybreaks
    \begin{align*}
\gamma (h_1) & = J_1^2 \, c (h_1) +  J_2^2 \, {c(h_1)}^2 \\
\gamma_3 (h_1, h_2) & =  \sqrt{2} J_1^2 J_2 
\, \left( c(h_1) c(h_2) + c(h_1) c(h_1 - h_2) + c(h_2) c(h_1 - h_2)
 \right) \\
 & +  {( \sqrt{2} J_2 )}^3 \, c(h_1) \, c(h_2) \, c(h_1 - h_2). 
\end{align*}
}
The expression for 
$\gamma_4 (h_1,h_2,h_3)$ is very long, and is given in the Supplement
\citet{GML2024}.

\subsection{Simulation Results}
For each of these three processes, we 
construct bispectral means 
from the following weight functions:
\begin{itemize}
    \item $g_1(\lambda_1, \lambda_2) = (4 \pi)^{-2}cos(3 \lambda_1)cos(\lambda_2)$ (as discussed in Example 1)
    \item $g_2(\lambda_1, \lambda_2) = \mathbb{I}_{[-0.2,0.2]} (\lambda_1)\mathbb{I}_{[-0.5,0.5]}(\lambda_2)$ (as discussed in Example 4)
    \item $g(\lambda_1, \lambda_2) = 1 - \sqrt{\frac{\lambda_1^2 + \lambda_2^2}{2}}$
\end{itemize}

\begin{figure}[htbp]
    \centering
    \includegraphics[width = 0.25\textwidth]{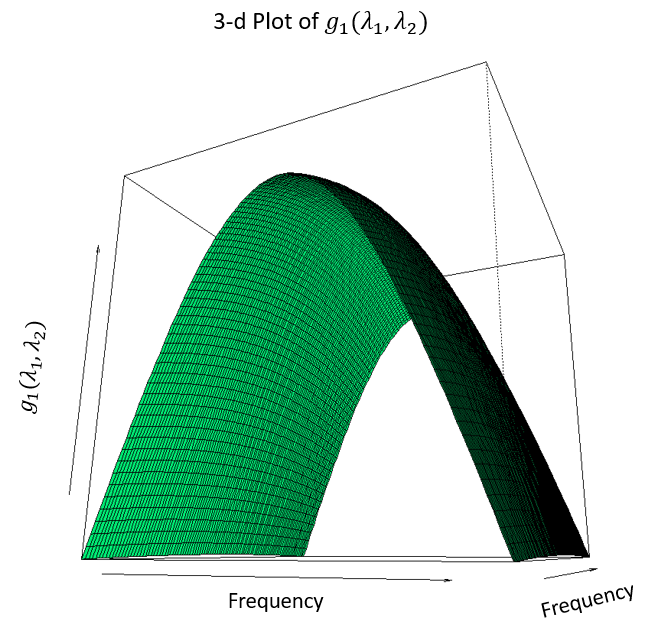}
    \includegraphics[width = 0.25\textwidth]{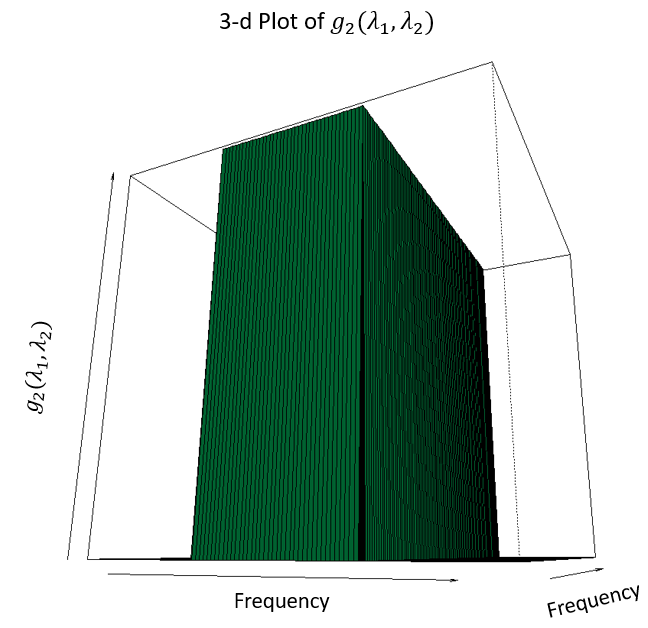}
    \includegraphics[width = 0.25\textwidth]{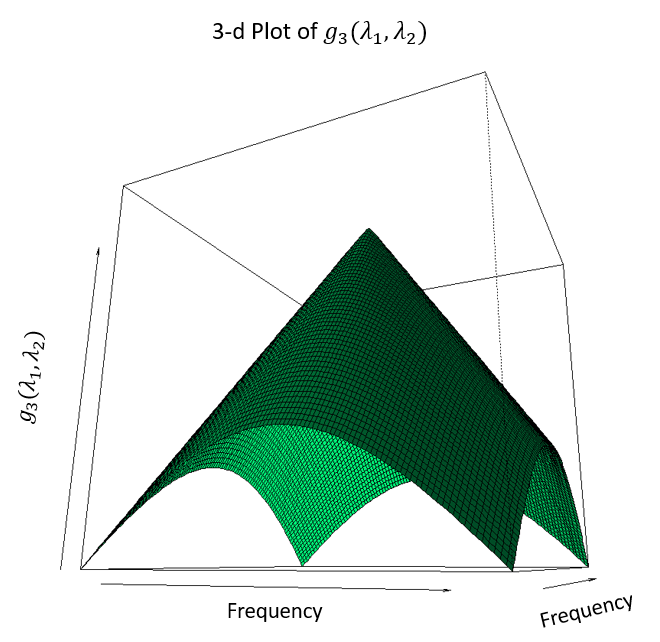}
    \caption{3-d plot of the weight functions given in the text. Polyspectral means with different weight functions   provide different features of the time series.}
    \label{fig:3dWeightFun}
\end{figure}

The 3-d plots of the weight functions are given in Figure \ref{fig:3dWeightFun}. 
Each choice of $g$ determines the
parameter $M_g (f_2)$,
and we also compute the true asymptotic variance $V$ as given in (\ref{eq:asymp-var}).
%, using the known weight functions and known polyspectra.
Next, we compute 
$\widehat{M_g (f_2)}$ from simulations
of length $T=100$, with $1000$ replications. 
The sample variance of these $1000$ values of $\widehat{M_g(f_2)}$ provides an estimate
$\widehat{V}$ of $V$; repeating the entire
procedure another $1000$ times yields
variance estimates  
%we will obtain an estimate of $V$ by taking the variance, say $\hat{V} = \frac{1}{1000}\sum_{i=1}^{1000} \left( \widehat{M_g(f_2)}_i - \overline{\widehat{M_g(f_2)}_i}\right)^2$ 
% Furthermore, we compute $\hat{V}$ 1000 times by the process described before, thereby getting estimates
% the bispectral mean estimate $\widehat{M_g (f_2)}$  for each simulation using (\ref{eq:polyHat}), and the variance is computed as an estimate 
$\hat{V}_i$ for $1 \leq i \leq 1000$.
 %of the proposed theoretical asymptotic variance. 
%This is replicated 1000 times, to get $\hat{V}_1, \ldots, \hat{V}_{1000}$.
The Mean Squared Error (MSE) and scaled MSE is then computed using the formulas $\sum_{i=1}^{1000} \left(\hat{V}_i - V \right)^2$ and $\sum_{i=1}^{1000} \left(\hat{V}_i/V - 1 \right)^2$, respectively.
The results of the simulation are given in Table \ref{table1}: for all cases, all models, and all weight functions, the computed asymptotic variance 
%in (\ref{eq:asymp-var}) 
is close  in terms of MSE to the 
simulated quantities.

\begin{table}[htbp]
\caption{Mean Squared Error (MSE) and scaled MSE of asymptotic variance of polyspectral mean for different processes and different weight functions (Sample Size is $100$).}
\label{table1}
\begin{center}
\begin{tabular}{|c|c|c|c|} 
\hline
Model & g($\underline{\lambda}$) & MSE & Scaled MSE \\
\hline
\multirow{3}{15em}{Model 1: AR(2) with $\mbox{Exp}(1) - 1$ errors} & $g_1(\lambda_1, \lambda_2)$  & 0.003 & 0.12 \\ 
& $g_2(\lambda_1, \lambda_2)$   & 0.023 & 0.19 \\ 
& $g_3(\lambda_1, \lambda_2)$    & 0.017 & 0.26 \\ 
\hline
% \multirow{3}{4em}{AR(2) with $\chi^2_2 -2$} & $g(\lambda) = 2 cos(2 \lambda)$  & &  \\ 
% & $g(\lambda)= \mathbb{I}_{[0,2]}(\lambda)$ & cell6 &\\ 
% & $g(\lambda) = \frac{4}{3}(1 - \lvert \lambda \rvert)/(2 \pi)$ & cell9 & \\ \hline

\multirow{3}{15em}{Model 2: AR(2) with $\chi^2_4 - 4$ errors} & $g_1(\lambda_1, \lambda_2)$  & $<$ 0.005 &  0.07   \\ 
& $g_2(\lambda_1, \lambda_2)$ & 0.011 & 0.15  \\ 
& $g_3(\lambda_1, \lambda_2)$ & 0.073 & 0.82  \\ \hline

\multirow{3}{15em}{Model 3: ARMA(2,1) with $\mbox{Exp}(1) - 1$ errors} & $g_1(\lambda_1, \lambda_2)$  & 0.023  & 0.39   \\ 
& $g_2(\lambda_1, \lambda_2)$ & 0.010 & 0.40  \\ 
& $g_3(\lambda_1, \lambda_2)$ & 0.027 & 0.14 \\ \hline

\multirow{3}{15em}{Model 4: ARMA(2,1) with $\chi^2_4 - 4$ errors} & $g_1(\lambda_1, \lambda_2)$  & 0.034 & 0.15   \\ 
& $g_2(\lambda_1, \lambda_2)$ & $<$ 0.005 & 0.28  \\ 
& $g_3(\lambda_1, \lambda_2)$ & 0.009 & 1.07  \\ \hline

\multirow{3}{15em}{Model 5: Squared Hermite process} & $g_1(\lambda_1, \lambda_2)$ & 0.081 & 0.27 \\
& $g_2(\lambda_1, \lambda_2)$ & 0.037 & 0.14  \\ 
& $g_3(\lambda_1, \lambda_2)$ & 0.002 & 0.91  \\ \hline
% \multirow{3}{4em}{Exponential Hermite} & $g(\lambda) = 2 cos(2 \lambda)$  & &  \\ 
% & $g(\lambda)= \mathbb{I}_{[0,2]}(\lambda)$ & cell6 & \\ 
% & $g(\lambda) = \frac{4}{3}(1 - \lvert \lambda \rvert)/(2 \pi)$ & cell9 & \\ \hline
\end{tabular} 
\end{center}
\end{table}

\subsection{Power Curve for Linearity Test}

\begin{figure}[htbp]
    \centering
    \includegraphics[width =0.5\textwidth]{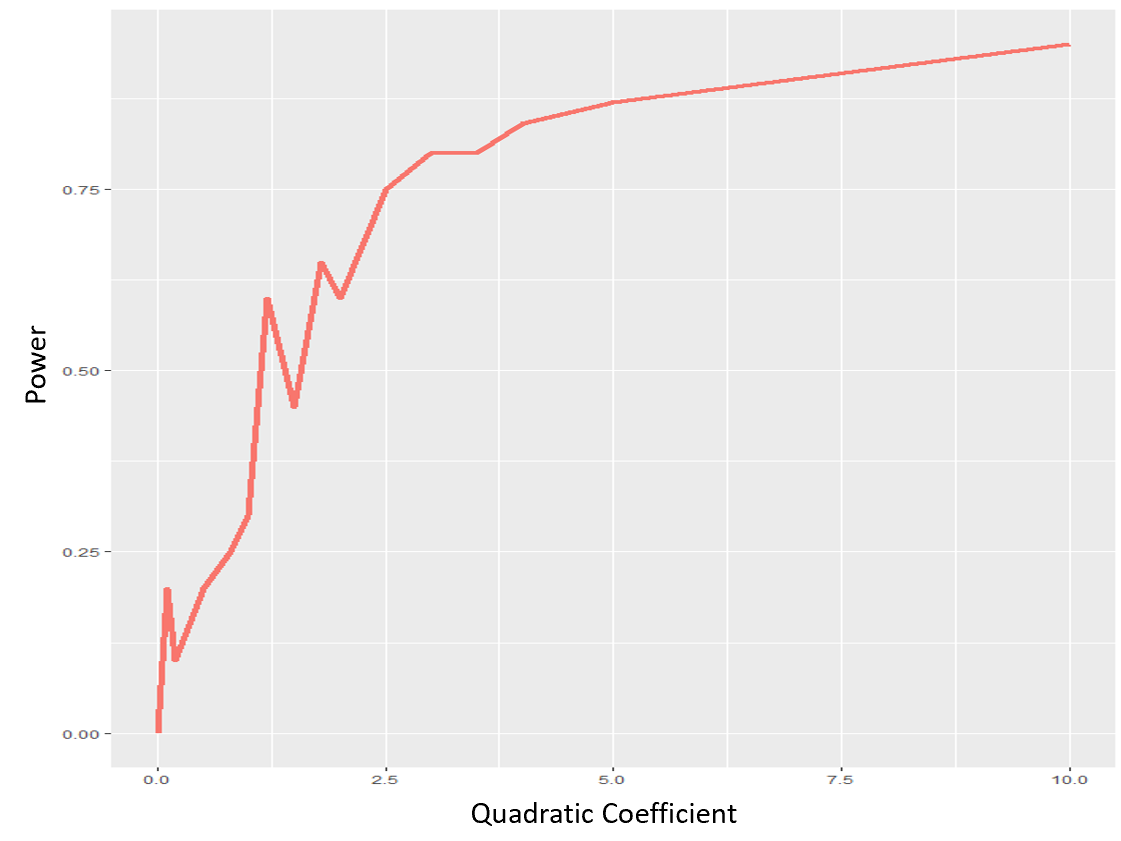}
    \caption{Power Curve of nonlinearity test for different choices of $\theta$ in the process
    (\ref{eq:quad-proc}).  Power is increasing 
    in $\theta $, which parameterizes deviation
    from linearity. }
    \label{fig:power}
\end{figure}

We   studied through simulations the power of our proposed linearity test (\ref{eq:blt-stat})
with $M=10$, generating the $\{ X_t \}$ process   from 
\begin{equation}
\label{eq:quad-proc}
     X_t = \epsilon_t + 0.4\epsilon_{t-1} + \theta \epsilon_{t-1}^2 - \theta, \quad \epsilon_t \sim \mathcal{N}(0, 1).
\end{equation}
Note that when $\theta =0$ this process
is linear, and hence corresponds to the null
hypothesis (\ref{eq:null-linear}) of our linearity test statistic.  Non-zero values of $\theta$ 
generate nonlinear effects, corresponding 
to the alternative hypothesis, and 
so we can determine power as a function of $\theta$.
 % Ideally,   $\theta$ will control the nonlinearity of the process, and as such increasing it should increase the power of the test, which can be seen in Figure \ref{fig:power}. 
 We generated $1000$ simulations of sample size $100$ of the process (\ref{eq:quad-proc}), letting $\theta$ range from $0$ to $10$. 
  The power curve (Figure \ref{fig:power}) 
 increases with $\theta$, as anticipated;
 there are a   few jags that are attributed to Monte Carlo error. 

% \begin{figure}
%     \centering
%     \includegraphics[width =0.7\textwidth]{table.png}
%     \caption{Caption}
%     \label{fig:my_label}
% \end{figure}

\section{Real Data Analysis}
\label{sec:data}

\subsection{Sunspot Data Linearity Test}

\begin{figure}[htbp]
    \centering
    \includegraphics[width = 0.5\textwidth]{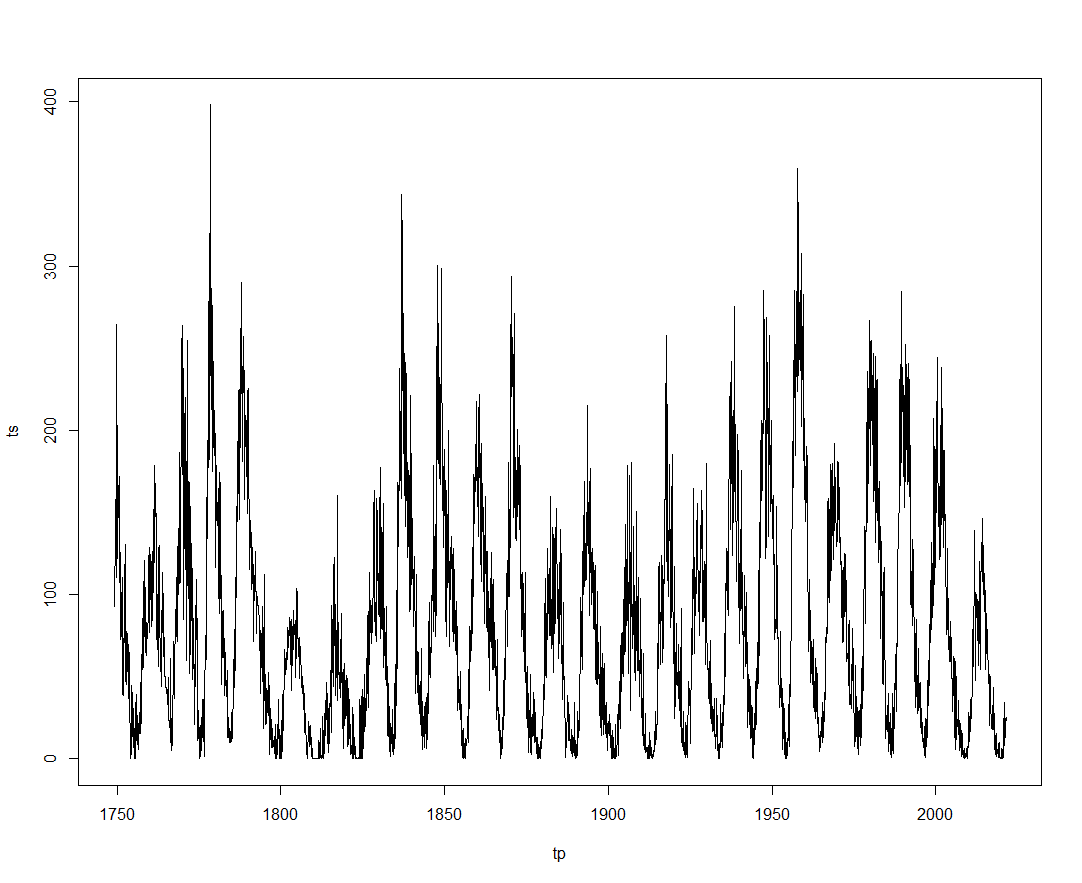}
    \caption{Time Series of Sunspot Data in a Monthly Interval (Source: \citet{sidc}).}
    \label{fig:sunspot}
\end{figure}

Sunspots are temporary phenomena on the Sun's photosphere that appear as spots relatively darker than surrounding areas. Those spots are regions of reduced surface temperature arising due to concentrations of magnetic field flux that inhibit convection. It is already known
that the solar activity follows a periodic pattern repeating every 11 years \cite{schwabe1844sonnenbeobachtungen}. Figure \ref{fig:sunspot} shows the time series of the sunspot data in a monthly interval \cite{sidc}.

\citet{ghosh2021} demonstrated that it is possible to gain significant improvement in prediction in sunspot data by using quadratic prediction instead of linear prediction. This means that there must be nonlinearity present   in the sunspot data. Here we will use the test statistic proposed in Section \ref{sec:application} to examine whether the sunspot data follows the particular linear process  studied in \citet{abdel2018statistical}, 
namely the order one autoregressive process given by
  $ X_t = 0.976 X_{t-1} + \epsilon_t$.  In the context
  of Section \ref{sec:application}, this means
  that $  \psi (z) = {(1 - .976 z)}^{-1}$
 and $\Psi (\lambda, \omega) =
 {(1 - .976 e^{-i \lambda} )}^{-1}
 {(1 - .976 e^{-i \omega})}^{-1}
 {(1 - .976 e^{i (\lambda + \omega)})}^{-1}$. 
Choosing $M=10$,  the test  statistic (\ref{eq:blt-stat})
    is $\mathcal{T}_{BLT} =  3740.057$,
which has a p-value $<0.005$ (Choosing $M=8$, the test statistic was $\mathcal{T}_{BLT} = 1945.213$ (p-value: $<0.005$) and choosing $M=5$, we get $\mathcal{T}_{BLT} = 1312.310$ (p-value: $<0.005$)). 
 Hence, our proposed test rejects the null hypothesis of linearity of the sunspot data, as assumed by \citet{abdel2018statistical}.
% A more general test is required to test whether the process is linear, which can be a topic for future research. 

\subsection{GDP Trend Clustering}

Time series clustering is a possible application
 of polyspectral means; by choosing various
 weighting functions, we can   extract different types of information from a time series. In this section, we have analyzed the Gross Domestic Product (GDP) of 136 countries (measured over the past 40 years), and have attempted to classify the growth trend of the GDP of different countries based on the bispectral means taken with different weight functions.  The right panel of Figure
 \ref{fig:gdp_figure} shows that there are
 different trend patterns for the various
 countries, and these diverse patterns are
 reflected in the differenced time series
 (left panel of Figure
 \ref{fig:gdp_figure}); the differencing
  ensures that the resulting growth rates are stationary.  
 % The time series of the GDP is first differenced
 %to ensure stationarity, and then scaled in order to extract just the trend information (shown in Figure \ref{fig:gdp_figure}).
 We investigate clustering of these time series
 through examining bispectral means of
 their growth rates. 
 % Hence, the clustering is based solely on how the GDP varied over time, and not the actual value of GDP. 
 The weight functions considered are:
%
% \begin{figure}[h!]
%     \centering
%     \includegraphics[width = 0.6\textwidth]{gdpplot.png}
%     \caption{Original and Differenced GDP trend for the countries over 40 years (1980-2020)}
%     \label{fig:gdp_figure}
% \end{figure}
%
\begin{itemize}
    \item $g(\lambda_1, \lambda_2) = \mathbb{I}( a<\lambda_1^2 + \lambda_2^2<b)$, where a and b are taken such that the interval $(0,1)$ is split into 10 segments. Hence, we have 10 such bispectral means for each annulus (as discussed in Example 4).
    \item $g(\lambda_1, \lambda_2) = (\pi- \lvert \lambda_1 \rvert)(\pi - \lvert \lambda_2 \rvert)$ (as discussed in Example 5).
    \item $g(\lambda_1, \lambda_2) = (2 \pi)^{-2} cos(3 \lambda_1)cos(\lambda_2)$ (as discussed in Example 1).
\end{itemize}
\begin{figure}[htbp]
    \centering
    \includegraphics[width = 0.6\textwidth]{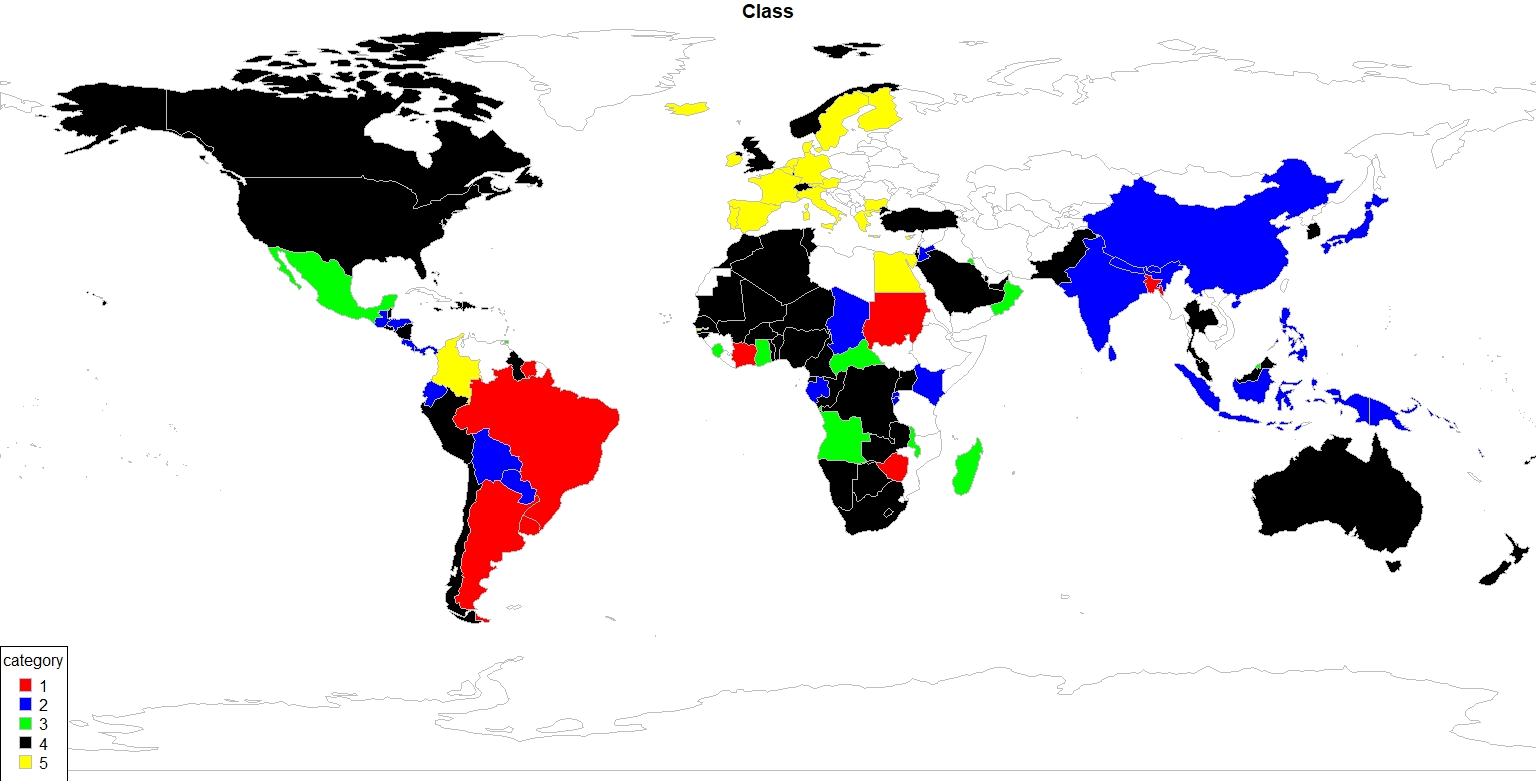}
    \caption{Clustering of 136 countries
    according to 12 polyspectral means 
    estimated from annual GDP growth rates.}
    \label{fig:cluster_world}
\end{figure}
The 12 computed bispectral means are then used to classify the 136 countries into 5 classes, as shown in Figure \ref{fig:cluster_world}, using k-means clustering. The features considered are all the bispectral means obtained using the given weight functions.  The developed EU countries all fall into the same category. The South Asian and South-East countries form another category, while all types of categories are found in Africa and South America. Some developed countries, like the USA and Canada,
are members of the same category as undeveloped
countries such as Benin, and this is likely due
to similar patterns of GDP growth. In other words, the classification algorithm classifies countries based on the growth trend of the GDP rather than the actual GDP, and the fact that USA and Benin are in the same group might be due to the fact that both countries have recently shown similar growth trend in GDP, the reasons behind which might be of practical interest, and can be explored in further studies.

\section{Proofs}
\label{sec:proof}

\subsection{Proof of Proposition \ref{prop1}}
\begin{proof}
\label{proof:prop1}
 By direct calculation, the expected value of
 $\widehat{M_g(f_k)})$ is
\[
    E\left(\widehat{M_g(f_k)})\right) = (2 \pi)^k T^{-k-1}\sum_{\tilde{\underline{\lambda}}} E\left( d\left(\tilde{\lambda_1}\right) \ldots d\left(\tilde{\lambda_k},  \right)d\left(-[\tilde{\underline{\lambda}}]  \right)\right)g\left(\tilde{\underline{\lambda}} \right). 
\]
Now, using results from \citet{brillinger1967asymptotic}
 we obtain 
\[
    E\left( d\left(\tilde{\lambda_1}\right) \ldots d\left(\tilde{\lambda_k}\right)d\left(-[\tilde{\underline{\lambda}}]  \right)\right)
    % g\left(\underline{\lambda} \right)
    = \sum_{\sigma \in \tau_n} \prod_{b \in \sigma} C_b.
\]
The term $C_b$ will be non-zero only when the sum of the Fourier frequencies inside the partition is $0$, which can only occur when the partition contains all the $\tilde{\lambda}$'s, since the Fourier Frequencies are assumed not to lie in any sub-manifold. Therefore, 
$ 
    E\left( d\left(\tilde{\lambda_1}\right) \ldots d\left(\tilde{\lambda_k} \right)d\left(-[\tilde{\underline{\lambda}}]  \right)\right)
    % g\left(\underline{\lambda} \right) 
 %   &= \sum_{\sigma \in \tau_n} \prod_{b \in \sigma} C_b \\
    = \Cum \left(d\left(\tilde{\lambda_1}\right),\ldots, d\left(\tilde{\lambda_k}\right), d\left(-[\tilde{\underline{\lambda}}] \right) \right) 
    =  T f\left(\tilde{\underline{\lambda}} \right) + O(1),
$ 
 where the last equality is a direct outcome of Lemma 1 of \citet{brillinger1967asymptotic}, and the error term is uniform for all $\{\tilde{\lambda_1}, \ldots,\tilde{\lambda_k}\}$.
Hence, using the convergence of a Riemann sum to 
an integral form (which is valid given the properties of the weight function), we obtain
$ 
    E\widehat{M_g(f)} = (2\pi)^k T^{-k-1} \sum_{\tilde{\underline{\lambda}}} \big[ T f(\tilde{\underline{\lambda}}) + O(1) \big] g(\tilde{\underline{\lambda}}) 
    = M_g\left(f\right) + o(1) $.
Note that by \eqref{assumptionA}, 
the $k$th order polyspectrum is bounded.
Hence,  it is easy to 
show that
the $o(1)$-term is indeed 
$ O\left(T^{-1}\right)$
when \eqref{eq:g-cond} holds. 
\end{proof}

\subsection{Proof of Proposition 2}
\label{proof:prop2}

\begin{proof}
%For the second cumulant, we can consider $\Cum(X,X), \Cum(X, \Bar{X})$ or $\Cum(\Bar{X}, \Bar{X})$. All the cases will have similar proofs, and we will only look at the circular cumulant (\citet{comon2010handbook}) $\Cum(X, \Bar{X})$ in this paper. We can write:
Recall that (\ref{cumulant1}) has been established,
and hence from (\ref{cumulant1}) and (\ref{eq:second}) we only need to examine
$\sum_{\sigma \in I_{2k+2}} \prod_{b \in \sigma} C_b$.
From Lemma 1 of \citet{brillinger1967asymptotic}, under the assumption that $\sum \lvert v_jc'_{a_1, \ldots, a_k}(v_1, \ldots, v_{k-1}) \rvert < \infty$, this 
term is  non-zero 
if and only if the sum of the frequencies is 0.  Recall that $\{\underline{\lambda}\}$ is a shorthand for  $\{\lambda_1, \ldots, \lambda_k, -[\underline{\lambda}] \}$.  Then it follows that
the terms $ \prod_{b \in \sigma} C_b$ will be non-zero in two possible partitioning of the vector $\{\{\underline\lambda\}, -\{\underline\omega\}\}$, as given in Section \ref{sec:results}. As we described in the referred section, there are two possible partitions, one with two sets containing elements only from $\{\underline{\lambda}\}$ and $-\{\underline{\omega}\}$, and others with a mixture of each vector, such that every set contains at least one element from each vector.

%  \begin{enumerate}
%     \item There are two partitions of $\{\{\underline{\lambda}\}, -\{\underline{\omega}\}\}$,  which is $\{\{\underline{\lambda}\}\}$ and $\{-\{\underline{\omega}\}\}$.
%    \item A second case would be where all the 
%    $\{\underline{\lambda}\}$ and 
%    $\{\underline{\omega}\}$ fall in the same partition. 
%      \item The other cases where the term would be non-zero are when the sums of all the elements within a block of a partition are $0$.
% \end{enumerate}     
    
In the first case, we have only one partition with two sets, namely, $\{\{\underline{\lambda}\}, -\{\underline{\omega}\}\}$. Let us denote this partition by $\sigma_0$. Since the sum of elements of both partitions is zero, both of them will contribute non-zero values to the product, rendering the product to be non-zero. Hence
$ \prod_{b \in \sigma_{0}} C_b$ can be written as
    \begin{align*}
        &= \Cum\left(d(\lambda_1), \ldots, d(\lambda_k), d\left(-[\underline{\lambda}] \right) \right) \Cum \left(d(-\omega_1), \ldots, d(-\omega_k), d\left([\underline{\omega}] \right) \right) \\
        &= \left( 
        % \left(2 \pi \right)^{k}
        Tf(\underline{\lambda})  + O(1) \right)\left( 
        % \left(2 \pi \right)^{k}
        Tf(-\underline{\omega})  + O(1) \right) \\
        &= 
        % \left(2 \pi \right)^{2k} 
        T^2 f(\underline{\lambda})f(-\underline{\omega}) + O(T).
    \end{align*} 
% \vspace{-1cm}
%
Hence, substituting the value of $ \prod_{b \in \sigma_{0}} C_b$ in the sum of (\ref{eq:second}), we obtain the following:
    \begin{align*}
        & T^{-2k-2}(2 \pi)^{2k} \sum_{ \tilde{\underline{\lambda}}, \tilde{\underline{\omega}}} g( \tilde{\underline{\lambda}} ) \overline{g(\tilde{\underline{\omega}})} \left\{ 
        % \left(2 \pi \right)^{2k}
        T^2 f( \tilde{\underline{\lambda}} ) f(-\tilde{\underline{\omega}}) + O(T) \right\}\\
        =& T^{-2k}\left(2 \pi \right)^{2k} \sum_{ \tilde{\underline{\lambda}}, \tilde{\underline{\omega}}} g( \tilde{\underline{\lambda}} ) \overline{g( \tilde{\underline{\omega}})} f(\tilde{\underline{\lambda}}) 
        f(-\tilde{\underline{\omega}}) + (2\pi)^{2k}T^{-2k}
        \sum_{\tilde{\underline{\lambda}}, \tilde{\underline{\omega}}} g(\tilde{\underline{\lambda}})
        \overline{g( \tilde{\underline{\omega}})}
        \cdot O\left(T^{-1}\right). 
    \end{align*}
Thus, the final expression is $M_g(f)\overline{M_g(f)} + O (T^{-1})$, which in view of  (\ref{cumulant1}) leaves a $O(T^{-1})$ term. The terms $O(T)$ and $O\left(T^{-1}\right)$ are uniform in $\{\lambda_1, \ldots, \lambda_k\}$, as given in \cite{brillinger1967asymptotic}. 

% \vspace{-2cm}
% \clearpage
% {\bf Need to say something about the uniformity of the $O(T)$
% and $O(T^{-1})$ terms. } 
The second case can have multiple possible partitions, where every set of every partition must contain a mixture of $\{\underline{\lambda}\}$ and $-\{\underline{\omega}\}$. Let us first consider the special case where we have only one set in the partition, i.e. the partition consists only of $\{\{\underline{\lambda}\}, -\{\underline{\omega}\}\}$. Let us call this partition as $\sigma_1$. Then, 
 \begin{align*}
    \prod_{b \in \sigma_{1}} C_b & =
         \Cum \left(d(\lambda_1), \ldots, d(\lambda_k), d\left(-[\underline{\lambda}]\right), d(-\omega_1), \ldots, d(-\omega_k), d\left([\underline{\omega}]\right) \right) \\
        &=  
        % \left(2 \pi \right)^{k}
        Tf_{2k+1}(\underline{\lambda}, -[\underline{\lambda}], -\underline{\omega})  + O(1). 
        % &= \left(2 \pi \right)^{2k} T^2 f(\underline{\lambda})f(-\underline{\omega}) + O(T)  
    \end{align*} 
    Hence, again substituting the value of $ \prod_{b \in \sigma_{1}} C_b$ in the sum (\ref{eq:second}) we obtain
    \begin{align*}
        & T^{-2k-2}(2 \pi)^{2k} \sum_{ \tilde{\underline{\lambda}}, \tilde{\underline{\omega}}} g( \tilde{\underline{\lambda}} ) \overline{g( \tilde{\underline{\omega}} )}
        \left\{
        % \left( 
        % \left(2 \pi \right)^{k}
        Tf_{2k+1}( \tilde{\underline{\lambda}},
        -[\tilde{\underline{\lambda}}],
        -\tilde{\underline{\omega}})
        + O(1) 
        % \right)
        \right\}\\
        &= T^{-2k-1}(2\pi)^{2k} \sum_{\tilde{\underline{\lambda}}, \tilde{\underline{\omega}}} g(\tilde{\underline{\lambda}}) \overline{g(\tilde{\underline{\omega}})} f_{2k+1} (\tilde{\underline{\lambda}},
        -[\tilde{\underline{\lambda}}],
        - \tilde{\underline{\omega}}) + (2\pi)^{2k}T^{-2k}
        \sum_{\tilde{\underline{\lambda}}, \tilde{\underline{\omega}}} g(\tilde{\underline{\lambda}}, \tilde{\underline{\omega}}) \cdot  O\left(T^{-2}\right) \\
        & = \Big[T^{-1}\int_{-\pi}^\pi \ldots \int_{-\pi}^\pi g(\underline{\lambda})\overline{g(\underline{\omega})}f_{2k+1}(\underline{\lambda},-[\underline{\lambda}],-\underline{\omega})d\underline{\lambda}d\underline{\omega}
       \Big] \big( 1+ o(1)\big) 
        + O(T^{-2}).
    \end{align*}
    % Thus, the final expression will be asymptotically:
    % \begin{equation}
    %     M_g(f_{2k+1}) + O\left( T^{-2}\right)
    % \end{equation}
Hence, the final expression will be a polyspectral mean of order $2k+1$, with a weight function depending on the original weight function. 
    
   The third case includes the situation where 
 if the blocks of a partition $\sigma$ are $b_1, b_2, \ldots, b_p$, and the $b_i$ block contains elements $\lambda_{b_{i1}}, \ldots, \lambda_{b_{il}}, -\omega_{b_{i1}}, \ldots, -\omega_{b_{im}}$, then the sums of all those elements must be $0$ for each of the blocks $b_i$. This is because if any of the blocks has elements which do not add up to $0$, the corresponding $C_b$ will be zero, and hence $\prod_{b \in \sigma} C_b$ will become zero for that partition. Now, we have ensured that no subset of the Fourier frequencies lies in a sub-manifold, i.e., for no subset of $\lambda_1, \ldots, \lambda_k$ can the sum   be $0$. Hence, for one partition to have a non-zero value, all of its blocks must contain at least one element from each of $\lambda$ and $\omega$.  The only other case is covered in the first case, where the partition contains all elements of $\lambda$ and $\omega$. Now that we know that every block of the partition contains a mixture of $\{\underline{\lambda}\}$ and $\{\underline{\omega}\}$, let us call each mixture a constraint since each mixture should have a sum 0, and hence would produce a constraint and thereby reducing the degrees of freedom by 1. Suppose we have $m$ such constraints, i.e., $m$ linear combinations of $\lambda_k$ and $\omega_k$ are $0$. This also means we have $m$ blocks. Then we will have a $T^{m}$ factor from the integrand since every block would contribute a $T$ from the cumulant, and a $O(1)$ from the residual. The sum will hence run for only the above cases, i.e., only when the linear combinations of the elements of partitions are equal to 0. In other words, only when there exists 
 $A, B \in \widetilde{\mathcal{B}}_{m,k+1}$  such that
 $A \{\underline{\lambda}\} = B \{\underline{\omega}\}$,
    where $\{\underline{\lambda}\} = \{\lambda_1, \ldots, \lambda_k, -[\underline{\lambda}]\}$. 
%The same argument applies for $\{ \underline{\omega}\}$.
Recall that the rows and columns of $A$ (and $B$)  have the properties delineated in the definition of $\widetilde{\mathcal{B}}_{m,k+1}$, arising from the properties of the partition.

% must be disjoint -- this means that the rows of A (and B) cannot have a 1 in the same column position -- and contain $1$ or $0$, and not all $0$. Also, every column must have exactly one 1. This is because: 
%     \begin{itemize}
%     \item A Fourier Frequency $\lambda$ can be only in one partition, and hence the rows of A must be disjoint
%     \item All columns must have at least one $1$, since a Fourier frequency $\lambda$ must be in at least one partition. 
%     \item Note that it follows that a column cannot have more than one $1$, since that would violate the disjoint row property. Hence, every column must have exactly one $1$.
%     \end{itemize}
%     An example is furnished by
% \begin{equation*}
% A_{m,k+1} = 
% \begin{pmatrix}
% 1 & 0 & \cdots & 1 \\
% 0 & 0 & \cdots & 0 \\
% \vdots  & \vdots  & \ddots & \vdots  \\
% 0 & 1 & \cdots & 0 
% \end{pmatrix}_{m \times k},
% \end{equation*}
% %
Let $f_\alpha(\underline{\lambda}_\alpha)$, 
$r_{Aj} $, $r_{Bj} $, $r_j$,
 $\lambda_{r_{A_j}}$,  
 and  $\omega_{r_{B_j}}$ be as defined in Section \ref{sec:results}. Then   each term inside of the sum of (\ref{eq:second}) is given by  (recall 
 the notation $\tilde{f}(\{\underline{\lambda}\}) = f(\underline{\lambda})$ introduced
 previously) 
\begin{equation*}
    T^{m} \prod_{j=1}^m \tilde{f}_{r_j }\left( \lambda_{r_{A_j}} , \omega_{r_{B_j}} \right) + O\left(T^{m-1} \right).
\end{equation*}
Now, for each of these quantities, we have $m-1$ constraints ($m$ of which arise from each row
of the matrices in $\widetilde{\mathcal{B}}_{m,k+1}$,
 lessened by one  because of a redundant constraint
 that the entire sum is $0$). Then,  for each   term
 we have
\begin{align*}
    &T^{-2k-2} \sum_{\tilde{\underline{\lambda}}, \tilde{\underline{\omega}}} g\left( \tilde{\underline{\lambda}} \right)\overline{g\left( \tilde{\underline{\omega}} \right)}\left\{T^{m} \prod_{j=1}^m \tilde{f}_{r_j }\left( \lambda_{r_{A_j}} , \omega_{r_{B_j}} \right) + O\left(T^{m-1} \right) \right\} \\
    &= T^{-1}T^{-2k+m-1}\sum_{\tilde{\underline{\lambda}}, \tilde{\underline{\omega}}} g\left( \tilde{\underline{\lambda}} \right)\overline{g\left( \tilde{\underline{\omega}} \right)}\prod_{j=1}^m \tilde{f}_{r_j }\left( \lambda_{r_{A_j}} , \omega_{r_{B_j}} \right) + O\left(T^{-2} \right).
\end{align*}
This goes to 0 at rate  $T^{-1}$. Recall the
definitions of $L_m$, $\underline{l}$ and $\zeta_{\underline{l}}$ from Section \ref{sec:results}.
%
% Let $L_m = \left\{l_1, \ldots, l_m \right\}$ be such that $\sum_{j=1}^m l_j = k+1$, $l_j > 0$. For example, for $k=2$ and $m=2$, $L_2 = \{(1,2), (2,1)\}$.  Further suppose that
% \begin{equation*}
%     \zeta_{L_m} = \Bigg\{A_{m \times (k+1)}\Bigg\vert A_{i\cdot}  \textrm{ has } l_i \textrm{ 1's and every column has exactly one 1}\Bigg\}.
% \end{equation*}
%
% \begin{equation*}
%     \Xi_{L_m, L_n} = \left\{ \underline{\lambda}, \underline{\omega} \Bigg\vert A \underline{\lambda} + B \underline{\omega} = 0 \textrm{ for some } A \in \zeta_{L_m}, B \in \zeta_{L_m'} \right\}  
% \end{equation*}
%
%Then for every choice of $m$, we will have a corresponding $L_m$, and for each elements $l_m \in L_m$, we will have corresponding sets $\zeta_{\underline{l}^{(m)}}$.
Then we will have the summation in (\ref{eq:second}) running over all possible choices of matrices $A$ and $B$ within the sets $\zeta_{\underline{l}}$ for all choices of $\underline{l} \in L_m$, such that $A\{\underline{\lambda}\} = B \{\underline{\omega}\}$. Then each of the summands will be of the form
%
% $$ V =  \sum_{m=1}^{k+1}\sum_{ \underline{l}_m , \underline{l}_m' \in L_m} \sum_{A \in \zeta_{l_m}, B\in \zeta_{l_m'}} 
%     % \tau_{\underline{l}_m}\tau_{\underline{l}_m'}
%     \underbrace{\int_{-\pi}^\pi \ldots \int_{-\pi}^\pi}_{A\{\underline{\lambda}\} - B \{\underline{\omega}\} = 0} g(\underline{\lambda})\overline{g(\underline{\omega})} \prod_{j=1}^m  \tilde{f}_{r_j } \left(\lambda_{r_{A_j}} , \omega_{r_{B_j}} \right)d\underline{\lambda}d\underline{\omega}$$
%
\[
T^{-2k-2} T^m \sum_{ \tilde{\underline{\lambda}}, \tilde{\underline{\omega}} \vert A\{\underline{\lambda}\} = B \{\underline{\omega}\} }
 g\left(\tilde{\underline{\lambda}} \right)
\overline{g\left( \tilde{\underline{\omega}} \right)}\prod_{j=1}^m \tilde{f}_{r_j }\left(\lambda_{r_{A_j}}, \omega_{r_{B_j}} \right),
\]
    % &=  T^{-1} T^{-2k + m -1}\sum_{\underline{\lambda},\underline{\omega} \vert A\{\underline{\lambda}\} - B \{\underline{\omega}\} = 0}g\left(\underline{\lambda}\right)\overline{g\left(\underline{\omega}\right)}\prod_{j=1}^m \tilde{f}_{r_j }\left(\lambda_{r_{A_j}} , \omega_{r_{B_j}} \right) 
%
% Now, let's see how many matrices are there in $\zeta_{L_m}$. This would be:
%
% \begin{align*}
%     \tau_{\underline{l_m}} &= \binom{k+1}{l_1}\binom{k+1-l_1}{l_2}\binom{k+1-l_1-l_2}{l_3}\ldots \binom{k+1-l_1-\ldots l_{m-1}}{l_m} \\
%     &= \frac{(k+1)!}{l_1! \ldots l_m!}
% \end{align*}
%
 and the powers of $T$ can be expressed
 as $T^{-1} T^{-2k + m -1}$.  It follows that the second cumulant is of the form 
\[
   T^{-1} \sum_{m=1}^{k+1}\sum_{\underline{l} \in L_m}     T^{-2k + m -1}\sum_{ \tilde{\underline{\lambda}},
   \tilde{\underline{\omega}} \vert A \{\underline{\lambda}\} = B \{\underline{\omega}\} }
    g\left( \tilde{\underline{\lambda}} \right)
\overline{g\left( \tilde{\underline{\omega}} \right)}\prod_{j=1}^m \tilde{f}_{r_j }\left(\lambda_{r_{A_j}} , \omega_{r_{B_j}} \right) + O(T^{-2}),
\]
from which the stated formula for $V$
 (\ref{eq:asymp-var}) is obtained. 
    % &\sim  \sum_{m=1}^{k+1}\sum_{ \underline{l}_m , \underline{l}_m' \in L_m} \sum_{A \in \zeta_{l_m}, B\in \zeta_{l_m'}} 
    % % \tau_{\underline{l}_m}\tau_{\underline{l}_m'}
    % \underbrace{\int_{-\pi}^\pi \ldots \int_{-\pi}^\pi}_{A\{\underline{\lambda}\} - B \{\underline{\omega}\} = 0} g(\underline{\lambda})\overline{g(\underline{\omega})} \prod_{j=1}^m  \tilde{f}_{r_j } \left(\lambda_{r_{A_j}} , \omega_{r_{B_j}} \right)d\underline{\lambda}d\underline{\omega}.
The case $m=1$ is the second case mentioned earlier, which would give a polyspectral mean of order $2k+1$.
This completes the proof. 
\end{proof}

\subsection{Proof of Corollary 1}
\label{proof:corr1}

\begin{proof}
The proof follows from the fact that
\begin{align*}
    &Cov(\widehat{M_{g_1}(f_k)}, \overline{\widehat{M_{g_1}(f_k)})} \\
    &= E \left\{ T^{-k} \sum_{\tilde{\underline{\lambda}}_k}\hat{f}_k(\tilde{\underline{\lambda}}_k)g_1(\tilde{\underline{\lambda}}_k) T^{-k-1}\sum_{\tilde{\underline{\lambda}}_{k+1}}\hat{f}_{k+1}(\tilde{\underline{\lambda}}_{k+1})\overline{g_2(\tilde{\underline{\lambda}}_{k+1})} \right\} - M_{g_1}(f_k)M_{g_2}(f_{k+1}) \\
    &= T^{-2k-3} (2 \pi)^{2k+1} \sum_{\tilde{\underline{\lambda}}_k, \tilde{\underline{\omega}}_{k+1}} g_1(\tilde{\underline{\lambda}}_k) \overline{g_2(\tilde{\underline{\lambda}}_{k+1})} \sum_{\sigma \in I_{2k+3}} \prod_{b \in \sigma} C_b  - M_{g_1}(f_k)M_{g_2}(f_{k+1}).
\end{align*}
% \begin{align*}
%     A &= T^{-2k-1} E \left\{ \sum_{\underline{\lambda}_k, \underline{\omega}_{k+1}} \hat{f}_k(\underline{\lambda}_k)\hat{f}_{k+1}(\underline{\omega}_{k+1}) g_1(\underline{\lambda}_k)g_2(\underline{\omega}_{k+1}) \right\} \\
%     &= T^{-2k - 3} (2 \pi)^{-2k -1} E \Bigg\{\sum_{\underline{\lambda}_k, \underline{\omega}_{k+1}} d(\lambda_1) \ldots d(\lambda_k)d\left( - [\underline{\lambda}]  \right). \\
%     &d(-\omega_1) \ldots d(  [\underline{\omega}] )d \left( \sum_{j=1}^k \omega_j \right) g_1(\underline{\lambda}_k)g_2(\underline{\omega}_{k+1}) \Bigg\} \\
%     &= T^{-2k-3} (2 \pi)^{-2k-1} \sum_{\underline{\lambda}_k, \underline{\omega}_{k+1}} g_1(\underline{\lambda}_k)g_2(\underline{\omega}_{k+1})  E \left\{ \prod_{\lambda \in \underline{\lambda}_k'} \prod_{\omega \in \underline{\omega}_{k+1}'}d(\underline{\lambda})d(\underline{\omega}) \right\} \\
%     &= T^{-2k-3} (2 \pi)^{-2k-1} \sum_{\underline{\lambda}_k, \underline{\omega}_{k+1}} g_1(\underline{\lambda}_k)g_2(\underline{\omega}_{k+1})  \sum_{\sigma \in I_{2k+3}} \prod_{b \in \sigma} Cum(b)
% \end{align*}
%
In this case, the proof would be similar to previous calculations, except the $A$ matrix will be of dimension $m \times (k+1)$ and the $B$ matrix will be of dimension $m \times (k+2)$. The proof can be easily extended to the covariance between polyspectral means of any two orders.
\end{proof}

\subsection{Proof of Theorem \ref{theo1}}
\label{proof:theorem1}

\begin{proof}
Propositions 1 and 2 give the mean and variance of the limiting distribution. The only thing remaining to prove asymptotic normality is to show that the higher-order cumulants of the scaled transformation go to zero as n goes to $\infty$. Recalling that $\{\underline{\lambda}\} = (\lambda_1, \ldots, \lambda_k, -[\underline{\lambda}])$,   we can write (letting $Cum_r$ denote the $r^{th}$ order joint cumulant)
\begin{align*}
    Cum_r(\widehat{M_g(f)}) &= Cum_r
    % \kappa
    (\underbrace{\widehat{M_g(f)}, \ldots, \widehat{M_g(f)}}_{\textrm{r times}}) \\
    &= Cum
    % \kappa
    \left(T^{-k-1} \sum_{\underline{\lambda}_1} g(\tilde{\underline{\lambda}_1})\prod d\{\tilde{\underline{\lambda_{1}}}\}, \ldots, T^{-k-1}\sum_{\tilde{\underline{\lambda}_r}} g(\tilde{\underline{\lambda}_r})\prod d\{\tilde{\underline{\lambda_{r}}}\} \right) \\
    &=  T^{-rk-r} \sum_{\tilde{\underline{\lambda_1}}}\ldots \sum_{\tilde{\underline{\lambda}_r}} g(\tilde{\underline{\lambda}_1})\ldots g(\tilde{\underline{\lambda}_r}) 
   \, 
    Cum \left( \prod d\{\tilde{\underline{\lambda_{1}}}\} , \ldots, \prod d\{\tilde{\underline{\lambda_{r}}}\} \right).
\end{align*}
%
% From Theorem 2.3.3 of \citet{brillinger2001time}, we know that
% \begin{align*}
%     \kappa\left(\prod_{j=1}^{J_1} X_{1j}, \ldots, \prod_{j=1}^{J_I} X_{Ij}\right) &= \sum_{\nu} \kappa(X_{ij}; ij \in \nu_1) \ldots \kappa(X_{ij} ; ij \in \mu_p), 
% \end{align*}
% where the summation is over all indecomposible partitions $\nu = \nu_1 \cup \ldots \cup \nu_p$ of the following table:
% \[
% \centering
% \begin{pmatrix}
% (1,1) & \cdots & (1,J_1) \\
% \cdot & & \cdot \\
% \cdot & & \cdot \\
% \cdot & & \cdot \\
% (I,1) & \cdots & (I,J_I) \\
% \end{pmatrix}
% \]
Applying  Theorem 2.3.3 of \citet{brillinger2001time},
we find that
\begin{align*}
    % \kappa
    Cum\left( \prod
    % _{j=1}^{k+1}
    d\{\tilde{\underline{\lambda_{1}}}\}_j , \ldots, \prod
    % _{j=1}^{k+1}
    d\{\tilde{\underline{\lambda_{r}}}\}_j \right) &= \sum_{\nu} 
    % \kappa
    Cum\left(d(\tilde{\lambda_{ij}}) : ij \in \nu_1\right) \ldots 
    % \kappa
    Cum\left(d(\tilde{\lambda_{ij}}) : ij \in \nu_p\right),
\end{align*}
where the summation is over all indecomposible partitions $\nu = \nu_1 \cup \ldots \cup \nu_p$ of the following table:
\[
\centering
\begin{pmatrix}
(1,1) & \cdots & (1,k+1) \\
\cdot & & \cdot \\
\cdot & & \cdot \\
\cdot & & \cdot \\
(r,1) & \cdots & (r,k+1) \\
\end{pmatrix}.
\]
From the expression of cumulants, we know that each of those 
% $\kappa$'s 
cumulants is non-zero only when the sum of the corresponding Fourier frequencies inside is equal to zero,
and can contribute terms of order at most $T$ for each of the partitions, i.e., the highest order from each of these elements is $T^p$. Hence, to obtain the leading term we would need to take the maximum number of such partitions such that the sum of $\lambda$'s in the partitions is zero. By our definition of $\widehat{M_g(f)}$, no subset of the $\underline{\lambda}_i$'s fall in a sub-manifold. Hence, each of these partitions must contain either all distinct $\{\underline{\lambda}\}s$, or a mixture of $\lambda$'s. However, if they were all distinct then the partition would not be indecomposable. Hence, we are left with only mixture of $\lambda$'s. Also, an indecomposable partition would mean that we don't have any single $\underline{\lambda}$, i.e. all the partitions must be a mixture of all the $\lambda$'s. This is encoded by the matrix notation   described earlier. Here we just have $r$ matrices $A_1, \ldots A_r$ such that $A_1 \underline{\lambda}_1 + \ldots + A_r \underline{\lambda}_r = 0_{p \times 1}$. Hence, $p$ partitions would give $p$ constraints, thereby giving an order of $T^p$, and with these come $p -1$ constraints. Therefore, the sum would need to be weighted by $T^{-rk + p - 1}$ in order to be non-negligible asymptotically, and the remainder term would be $T^{-r+1}$. Thus, the cumulants of order higher than $r$ will go to $0$ at a rate faster than $T^{-1}$, thereby proving the asymptotic normality of the scaled transformation.

\end{proof}

\subsection{Proof of Theorem \ref{thm:lin-test}}
\label{theorem2Proof}

\begin{proof}
Recall  that $\mathcal{T}_{BLT}$  is defined by (\ref{eq:blt-stat}), and  $g_{j,k} = (x_1, x_2) = \exp \{ \iota j x_1 + \iota k x_2 \}/\Psi(x_1, x_2)$. 
% the test statistic after scaling by $V_{j,k}$ (the asymptotic variance $V$ given by (\ref{eq:asymp-var}))
% can be written
% \[
%     \mathcal{T}_{BLT} = \sum_{(j,k) \neq (0,0),0 \leq j,k \leq M}  \frac{T \lvert \widehat{ M_{g_{j,k}}(f)}\rvert^2}{V_{j,k}}.
% \]
%
First we need to show that  $M_{g_{j,k}}(f)$ is $0$;
since under the null hypothesis $f(\underline{\lambda})$ is a constant multiple of the denominator of $g_{j,k}(\underline{x})$,
 and $(j,k) \neq (0,0)$, the integral is zero.  Hence,  the asymptotic distribution of $\sqrt{ T}   \widehat{M_{g_{j,k}}(f)}  / \sqrt{V_{g_{j,k},g_{j,k}} }$ under $H_0$ is $\mathcal{N}(0,1)$, using Theorem \ref{theo1p} (with $k=2$).

Let us next consider the joint distribution of $\left(\widehat{M_{g_{j_1, k_1}}(f)}, \widehat{M_{g_{j_2, k_2}}(f)} \right)$. Suppose we want to find the limiting distribution of $\upsilon = a_1 \widehat{M_{g_{j_1, k_1}}(f)} + a_2 \widehat{M_{g_{j_2, k_2}}(f)}$
 for any scalars $a_1$ and $a_2$. We can write:
\begin{align*}
    \upsilon &= a_1 (2\pi)^{-k}T^{-k} \sum_{\tilde{\underline{\lambda}}^{(a_1)} }
 T^{-1}d\left(\tilde{\lambda}^{(a_1)}_1 \right)\ldots d\left(\tilde{\lambda}^{(a_1)}_k \right)d\left(-\sum_{j=1}^k \tilde{\lambda}^{(a_1)}_j \right) g_{j_1, k_1}\left(\tilde{\underline{\lambda}}^{(a_1)} \right) \\
    & \quad  + a_2 (2\pi)^{-k}T^{-k} \sum_{\tilde{\underline{\lambda}}^{(a_2)}} T^{-1}d\left(\tilde{\lambda}^{(a_2)}_1 \right)\ldots d\left(\tilde{\lambda}^{(a_2)}_k \right)d\left(-\sum_{j=1}^k \tilde{\lambda}^{(a_2)}_j \right)g_{j_2, k_2}\left(\tilde{\underline{\lambda}}^{(a_2)}\right) \\
    &= (2 \pi)^{-k}T^{-k}\sum_{\tilde{\underline{\lambda}}} T^{-1} d\left(\tilde{\lambda_1} \right)\ldots d\left(\tilde{\lambda_k} \right)d\left(-\sum_{j=1}^k \tilde{\lambda_j}\right)h(\tilde{\underline{\lambda}}),
\end{align*}
where $h(\underline{\lambda}) = a_1 g_{j_1, k_1}(\underline{\lambda}) + a_2 g_{j_1,k_1}(\underline{\lambda})$.  Since we have already established the asymptotic normality of polyspectral mean for any function, it follows that any linear combination is asymptotically normal. Therefore, the joint distribution is asymptotically normal.

Hence,   we now need to find the asymptotic covariance between $\widehat{M_{g_{j_1, k_1}}(f)}$ and $\widehat{M_{g_{j_2, k_2}}(f)}$. 
%
% \begin{align*}
%     E \left\{ \widehat{M_{g_{j_1, k_1}}(f)} \overline{\widehat{M_{g_{j_2, k_2}}(f)}} \right\} - M_{g_{j_1,k_1}}(f)M_{g_{j_2,k_2}}(f)
% \end{align*}
%
As we saw earlier, one term that would arise
corresponds to the partition consisting of
the sets $\underline{\lambda}$ and $\underline{\omega}$.  We only need to show that there are no other partitions  in this case. The covariance will be similar to the form of $V$ established in Theorem \ref{theo1}, as 
given in (\ref{covTerm}) of
Corollary \ref{covCorr}.
% %
% \begin{align}
% % \label{covTerm}
%      Cov\left(\widehat{M_{g_1}(f_k)}, \widehat{M_{g_2}(f_k)}\right) &=  \sum_{m=1}^{k+1}\sum_{ \underline{l}_m , \underline{l}_m' \in L_m} \sum_{A \in \zeta_{l_m}, B\in \zeta_{l_m'}} 
%     % \tau_{\underline{l}_m}\tau_{\underline{l}_m'}
%     \underbrace{\int_{-\pi}^\pi \ldots \int_{-\pi}^\pi}_{A\{\underline{\lambda}\} - B \{\underline{\omega}\} = 0} g_1(\underline{\lambda})\overline{g_2(\underline{\omega})} \\
%     &\prod_{j=1}^m  \tilde{f}_{r_j } \left(\lambda_{r_{A_j}} , \omega_{r_{B_j}} \right)d\underline{\lambda}d\underline{\omega}
% \end{align}
We have finally established the joint distribution of $(M_{g_{j_1,k_1}},M_{g_{j_2,k_2}})$ as a multivariate normal with covariance matrix given by
\[
\begin{pmatrix}
V_{ g_{j_1,k_1}, g_{j_1,k_1} }
 &  V_{ g_{j_1,k_1}, g_{j_2,k_2} } \\
V_{ g_{j_2,k_2}, g_{j_1,k_1} }  & 
V_{ g_{j_2,k_2}, g_{j_2,k_2} } 
\end{pmatrix}
\] 
for any arbitrary $j_1, j_2, k_1, k_2$, where the matrix entries are defined via (\ref{covTerm}),
i.e., $V_{ g_{j_1,k_1}, g_{j_2,k_2} }$
is the covariance   between the bispectral means with weight functions $g_{j_1,k_1}$ and $g_{j_2,k_2}$. 
%Therefore $\upsilon$ is a multivariate normal distribution with a covariance matrix, say $V_\Upsilon$.

Next, consider the joint distribution of 
the $M^2-1$ vector $\sqrt{ T}   \widehat{M_{g_{j,k}}(f)}  / \sqrt{V_{g_{j,k},g_{j,k}} }$ such that
 $(j,k) \in \mathcal{v}$.
Then  with $\mathbb{CV}_{BLT}^{(M)}$
defined by (\ref{cvblt}),
the asymptotic distribution of   $\mathcal{T}_{BLT}$ under the null hypothesis can be expressed as $X^{\prime} X$, where $X \sim \mathcal{N}_{M^2-1}(0, \mathbb{CV}_{BLT}^{(M)})$.
%
%  Suppose $Y = C_{BLT}^{-\frac{1}{2}} X$. Then, $X^TAX = Y^T C_{BLT}^{\frac{1}{2}}A C_{BLT}^{\frac{1}{2}} Y$.
The covariance matrix has an eigenvalue decomposition
of the form $\mathbb{CV}_{BLT}^{(M)} = P^{\prime} \Lambda P$, where 
$\Lambda = \mbox{diag}(\nu_1, \ldots, \nu_{M^2-1})$ 
is a diagonal matrix of the eigenvalues of $\mathbb{CV}_{BLT}^{(M)}$.
Let $U = \Lambda^{-1/2} P X$, so that $ U \sim \mathcal{N}_{M^2-1}(0, I)$,
 where $I$ is the $M^2-1$-dimensional identity matrix.
Then, $X^{\prime} X = U^{\prime} \Lambda U = \sum_{j=1}^{M^2-1} \nu_j U_j^2$. Hence, the asymptotic distribution under the null hypothesis is
expressed as a weighted sum of $\chi^2_1$ random variables, where the weights are given by the eigenvalues of the covariance matrix $\mathbb{CV}_{BLT}^{(M)}$.
\end{proof}

\section*{Acknowledgements} This report is released to inform interested parties of research and to encourage discussion.  The views expressed on statistical issues are those of the authors and not those of the U.S. Census Bureau.  This work partially supported by NSF grants  DMS 1811998,
DMS 2131233,  DMS 2235457.

\bibliography{ref}
\bibliographystyle{rusnat}

\end{document}